\documentclass[reqno,11pt]{amsart}
\usepackage{amsmath,amsthm}
\usepackage{latexsym}
\usepackage{amssymb}
\usepackage[all,cmtip]{xy}
\newtheorem{Proposition}{Proposition}[section]
\newtheorem{Definition}[Proposition]{Definition}
\newtheorem{Lemma}[Proposition]{Lemma}
\newtheorem{Theorem}[Proposition]{Theorem}
\newtheorem{MainTheorem}{Theorem}

\newtheorem{Remark}[Proposition]{Remark}
\newtheorem{Example}[Proposition]{Example}
\DeclareMathOperator{\p}{\mathbb{P}}
\newcommand{\R}{\mathbb{R}}

\DeclareMathOperator{\cl}{cl}
\DeclareMathOperator{\spt}{spt}

\DeclareMathOperator{\GT}{GT}
\DeclareMathOperator{\WF}{WF}
\DeclareMathOperator{\Val}{Val}

\title{The product on smooth and generalized valuations}
\author{Semyon Alesker}
\author{Andreas Bernig}

\email{semyon@post.tau.ac.il}
\email{bernig@math.uni-frankfurt.de}

\date{\today}
\address{School of Mathematical Sciences, Tel Aviv University, 69978 Tel Aviv, Israel}
\address{Institut f\"ur Mathematik, Goethe-Universit\"at Frankfurt, Robert-Mayer-Str. 10, 60054
Frankfurt, Germany}
\begin{document}

\begin{abstract}
The product of smooth valuations on manifolds is described in terms
of differential forms, Gelfand transforms and blow-up spaces. It is
shown that the product extends partially to generalized valuations
and corresponds geometrically to transversal intersections. This
result is used to prove a general kinematic formula on compact rank
one symmetric spaces.
\end{abstract}

\thanks{ 
\\S.A. was partially supported by ISF grants 1369/04 and 701/08.\\ A.B. was supported
 by the Schweizerischer Nationalfonds grants PP002-114715/1 and SNF 200020-121506/1.}
\maketitle

\tableofcontents

\section{Introduction}

Roughly speaking, a {\it valuation} is a finitely-additive
functional on some system of sets. Classical examples are valuations
on convex sets or on convex polytopes. Recently, a detailed study of
valuations on manifolds was carried out in
\cite{ale05a,ale05b,ale05c,ale05d,ale06}. On a smooth manifold $X$,
one considers the system $\mathcal{P}$ of all compact submanifolds
with corners and calls a finitely additive functional a {\it
valuation} on $X$. Under an important, but technical assumption of
{\it smoothness}, there is a surprisingly rich algebraic structure
on the space of valuations on $X$. By \cite{ale05c,ale05d}, the
space $\mathcal{V}^\infty(X)$ of smooth valuations on $X$ is an
algebra satisfying a version of Poincar\'e duality.

Let us recall the main steps of that product construction. The
product of polynomial convex valuations is constructed in
\cite{ale04}, using the solution of P. McMullen's conjecture
\cite{ale01}. This product is then extended in \cite{ale05a} to
smooth valuations on affine space. Then the product of smooth
valuations on a general smooth manifold is obtained in \cite{ale05c}
by gluing together the products from \cite{ale05a} on affine spaces
in each coordinate chart. The major problem in this approach is to
show that the product is independent of various choices involved:
the construction of the product on affine spaces in \cite{ale05a}
uses a very special and non-unique way of presentation of smooth
valuations; the general case in \cite{ale05c} uses a choice of
affine coordinate charts. In both steps analytical difficulties have
arose, and in the second step \cite{ale05c} they were resolved using
some geometric measure theory.

In this paper we give a new construction of this product in
different terms. This construction is motivated by the previous one,
but has a number of advantages over it. First, it is more invariant,
in particular independent of the coordinate charts. Second, the
analytic difficulties arising in working with this construction seem
to be easier. This allows us to obtain some applications of the new
construction which will be discussed below.

Let us describe our main results in greater detail. Let $n:=\dim X$
and set $\p_X:=\p_+(T^*X)$ for the cosphere bundle of $X$. First, we
construct an exterior product, which associates to a smooth
valuation on a manifold $X$ and a smooth valuation on a manifold $Y$
some (non-smooth) valuation on $X \times Y$ and extends the exterior
product constructed in \cite{ale05a}. For the notion of transversal
sets and transversal valuations (denoted by $\mathcal{V}_{tr}(X
\times Y)$), we refer to Section \ref{sec_transversal}.

\begin{MainTheorem} \label{mthm_ext_prod}
Let $X$ and $Y$ be smooth manifolds. There exists a unique bilinear
map
\begin{displaymath}
\boxtimes:\mathcal{V}^\infty(X) \times \mathcal{V}^\infty(Y) \to \mathcal{V}_{tr}(X \times Y)
\end{displaymath}
with the following properties:
\begin{enumerate}
\item (Invariance under open embeddings) \label{item_invariance}
If $\phi:X' \to X$ and $\psi:Y' \to Y$ are smooth maps with open
images and which are diffeomorphisms onto their images, then for
$\mu_1 \in \mathcal{V}^\infty(X), \mu_2 \in \mathcal{V}^\infty(Y)$
and transversal $P \in \mathcal{P}(X' \times Y')$ we have
\begin{equation}
(\mu_1 \boxtimes \mu_2)((\phi \times \psi)(P))=(\phi^* \mu_1 \boxtimes \psi^* \mu_2)(P).
\end{equation}
\item (Affine case) \label{item_affine_case}
If $X$ and $Y$ are affine spaces, then $\mu_1 \boxtimes \mu_2$ is the exterior product introduced in \cite{ale05a}.
\end{enumerate}
\end{MainTheorem}

The product of two valuations on $X$ is the restriction of the
exterior product on the diagonal in $X \times X$. Our main theorem
describes the product explicitly in terms of differential forms and
Gelfand transforms. Below we assume for simplicity  of notation that
$X$ is oriented.

A smooth valuation $\mu$ on $X$ can be represented by a pair of smooth differential forms $(\omega,\phi) \in \Omega^{n-1}(\p_X) \times \Omega^n(X)$ as
\begin{displaymath}
 \mu(P)=\int_{N(P)} \omega+\int_P \phi
\end{displaymath}
for all compact manifolds with corners $P$. Here $N(P)$ is the conormal cycle of $P$ (which is an oriented closed Lipschitz manifold).

These forms are not unique. The pairs $(\omega,\phi)$ which induce the trivial valuation are characterized in \cite{bebr07} using a second-order differential operator on contact manifolds, which is called Rumin operator and denoted by $D$ (compare Theorem \ref{thm_bebr} in Section \ref{section_background}).

In Section \ref{sec_blow_up} we will construct a double fibration
\begin{displaymath}
\xymatrix{
& \bar \p \ar[dl]_{\bar p} \ar[dr]^{\bar \Phi} &  \\
\p_X & & \p_X \times_X \p_X}
\end{displaymath}
where $\bar \p$ is some blow-up space over $\p_X \times_X \p_X$ and $\bar p,\bar \Phi$ are natural projection maps.

This double fibration induces a Gelfand transform $\GT:
\Omega^*(\p_X \times_X \p_X) \to \Omega^{*-n}(\p_X)$ defined by
$\GT=\bar p_*\bar \Phi^*$, where $\bar\Phi^*$ is the pull-back
operation, and $p_*$ is the operation of integration along the
fibers. With $q_1,q_2:\p_X \times_X \p_X \to \p_X$ denoting the
canonical projections; $s:\p_X \to \p_X$ denoting the natural
involution map; and $\pi:\p_X \to X$ the projection map, our main
theorem is the following.

\begin{MainTheorem} \label{mthm_prod}
Let $X$ be an $n$-dimensional oriented manifold; let $\mu_i \in
\mathcal{V}^\infty(X)$ be represented by $(\omega_i,\phi_i) \in
\Omega^{n-1}(\p_X) \times \Omega^n(X)$. Then the product $\mu_1
\cdot \mu_2$ is represented by
\begin{align*}
\omega & = \GT(q_1^* \omega_1 \wedge q_2^* D\omega_2) + \omega_1 \wedge \pi^*\pi_* \omega_2 \in \Omega^{n-1}(\p_X), \nonumber\\
\phi & = \pi_* (\omega_1 \wedge s^* (D\omega_2+\pi^* \phi_2))+ \phi_1 \wedge \pi_* \omega_2 \in \Omega^n(X).
\end{align*}
\end{MainTheorem}

Notice that the assumption on orientation of $X$ may be removed, but
the notation will become more complicated. Using the new
construction of the product, we prove the following theorem.
\begin{MainTheorem} \label{mthm_functional_calc}
Let $\mu$ be a smooth valuation on an $n$-dimensional manifold $X$.
For an entire function $f(z)=\sum_{k=0}^\infty a_k z^k, a_k \in \mathbb{C}$, the power series
\begin{displaymath}
\sum_{k=0}^\infty a_k \mu^k
\end{displaymath}
converges to a smooth valuation denoted by $f(\mu)$.
\end{MainTheorem}

Since the product on valuations satisfies Poincar\'e duality, it is
possible to introduce a class $\mathcal{V}^{-\infty}(X)$ of {\it
generalized valuations} as a completion of $\mathcal{V}^\infty$ with
respect to some weak topology (see \cite{ale05d} for the details or
Section \ref{Ss:van-on-mflds} below). Every smooth valuation is also
a generalized valuation. Moreover, a manifold with corners $A$
induces a generalized valuation $\Xi_\mathcal{P}(A)$. It was
conjectured in \cite{ale06} that the product of smooth valuations
can be partially extended to generalized valuations and that this
extension corresponds to transversal intersections. We show a
version of this conjecture by proving the following two theorems.

\begin{MainTheorem}
 There exists a partial product on $\mathcal{V}^{-\infty}(X)$ which is commutative and associative and extends the product of smooth valuations.
\end{MainTheorem}

We refer to Theorem \ref{thm_prod_gen} for the precise conditions under which the product of two generalized valuations exists.

\begin{MainTheorem} \label{mthm_mflds_corners}
If $P^{(1)},P^{(2)}$ are compact submanifolds with corners which intersect transversally,
then the product of $\Xi_\mathcal{P}(P^{(1)})$ and $\Xi_\mathcal{P}(P^{(2)})$ exists and equals $\Xi_\mathcal{P}(P^{(1)} \cap P^{(2)})$.
\end{MainTheorem}

In the case where $X$ is real-analytic and $P^{(1)},P^{(2)}$ are subanalytic,
we conjecture that an analogous statement holds true, but we were not able to prove it.

The main application of this theorem is to the integral geometry of compact rank one symmetric spaces. If $X=G/H$ is a compact rank one symmetric space, then the space $\mathcal{V}^\infty(X)^G$ of smooth $G$-invariant
valuations has a finite basis $\varphi_1,\ldots,\varphi_N$ and there
are kinematic formulas
\begin{displaymath}
 \int_G \mu(P_1 \cap gP_2)dg=\sum_{i,j=1}^N c_{ij}^\mu \varphi_i(P_1)\varphi_j(P_2), \quad P_1,P_2 \in \mathcal{P}.
\end{displaymath}

They can be encoded by the map
\begin{align*}
 k_G: \mathcal{V}^\infty(X)^G & \to \mathcal{V}^\infty(X)^G \otimes \mathcal{V}^\infty(X)^G\\
\mu & \mapsto \sum_{i,j=1}^N c_{ij}^\mu \varphi_i \otimes \varphi_j.
\end{align*}

Let $m_G:\mathcal{V}^\infty(X)^{G} \otimes \mathcal{V}^\infty(X)^{G}
\to \mathcal{V}^\infty(X)^{G}$ denote the multiplication map.
Consider the bilinear map $\mathcal{V}^\infty(X)^{G} \otimes
\mathcal{V}^\infty(X)^{G} \to \mathbb{C}$ given by
$(\phi,\psi)\mapsto (\phi\cdot\psi)(X)$. By the Poincar\'e duality
this is a non-degenerate pairing. Consider the induced map
$p_G\colon\mathcal{V}^\infty(X)^{G} \to \mathcal{V}^\infty(X)^{G*}$.

\begin{MainTheorem}  \label{mthm_kinematic}
The following diagram commutes
\begin{displaymath}
\xymatrixcolsep{3pc}
\xymatrix{\mathcal{V}^\infty(X)^G \ar[d]^{k_G} \ar[r]^{p_G}& \mathcal{V}^\infty(X)^{G*} \ar[d]^{m_G^*}\\
\mathcal{V}^\infty(X)^G \otimes \mathcal{V}^\infty(X)^G \ar[r]^-{p_G \otimes p_G} & \mathcal{V}^\infty(X)^{G*} \otimes \mathcal{V}^\infty(X)^{G*}},
\end{displaymath}
\end{MainTheorem}

Thus the coefficients of the kinematic formulas
are closely related to the product structure on the space of smooth valuations.
A similar theorem in the affine, translation invariant case was proved in \cite{befu06}.
It was the starting point to  write down in an explicit form all kinematic formulas for the groups
$U(n)$ \cite{befu08}, $SU(n)$ \cite{be08a}, $G_2$ and $Spin(7)$ \cite{be08b}.

\begin{Remark}
Compact rank one symmetric spaces have been classified by E. Cartan (see e.g.
\cite{helgason}, \cite{bes87}). Cartan's list of simply connected compact rank one symmetric spaces is as follows:
\begin{displaymath}
S^n, \mathbb{CP}^n, \mathbb{HP}^n, \mathbb{C}a\mathbb{P}^2=F_4/Spin(9).
\end{displaymath}
\end{Remark}

\subsection*{Plan of the paper}

In Section \ref{section_background} we collect some important material on valuations and geometric measure theory. The blow-up spaces which will be needed in the paper are constructed in Section \ref{sec_blow_up}. In Section \ref{sec_transversal} we introduce transversal sets and transversal valuations and prove a formula for the Minkowski sum of a transversal set with a product set. The exterior product of smooth valuations is constructed in Section \ref{section_ext_prod}. The construction of the product is contained in Section \ref{section_prod}. Section \ref{sec_functional_calc} contains the proof of Theorem \ref{mthm_functional_calc}. The partial extension of the product to generalized valuations is contained in Section \ref{section_prod_gen}, where we give precise conditions under which the product exists. It is shown in Section \ref{sec_prod_man_corners} that these conditions are satisfied in the case of a transversal intersection of manifolds with corners and that Theorem \ref{mthm_mflds_corners} holds. The integral geometry of compact rank one symmetric  spaces is studied in Section \ref{section_kin_form}, which contains the proof of Theorem \ref{mthm_kinematic}.

\subsection*{Acknowledgments}
We wish to thank J. Fu for many interesting discussions related to
this work. The first author thanks J. Bernstein for introducing him
to the theory of wave front sets, and D. Kazhdan for very useful
stimulating  discussions. Part of the results in this paper were
proved during our mutual visits at the University of Fribourg and
the Institute of Advanced Studies in Jerusalem and we thank these
institutions for their hospitality. We thank F. Schuster for his
kind invitation for both of us to the Technical University of Vienna,
where we could finish the typing of
this manuscript.

We thank the anonymous referee for various very useful remarks on this text.

\section{Background}
\label{section_background}

\subsection{Smooth valuations on manifolds}\label{Ss:van-on-mflds}
The general reference for smooth valuations on manifolds is the series of papers \cite{ale05a,ale05b,ale05c,ale05d} as well as the survey \cite{ale06}.

A valuation on a smooth manifold $X$ is a functional
$\mu$ on the space $\mathcal{P}(X)$ of compact manifolds with corners which is finitely additive.
Roughly speaking, $\mu$ is called finitely additive if it satisfies an inclusion-exclusion principle whenever it makes sense,
see \cite{ale05b} for the precise definition.

Let us now suppose that $X$ is an oriented, $n$-dimensional manifold.
Every $P \in \mathcal{P}(X)$ admits a conormal cycle $N(P)$,
which is a closed Legendrian $n-1$-dimensional Lipschitz submanifold of the cosphere bundle $\p_X:=\p_+(T^*X)$. The orientation of $N(P)$ is fixed in such a way that
\begin{equation} \label{eq_image_down_normal_cycle}
 \pi_* N(P)=\partial P,
\end{equation}
where $\pi:\p_X \to X$ is the projection map.

A valuation $\mu$ on $X$ is called smooth if it can be represented by
differential forms $\omega \in \Omega^{n-1}(\p_X)$ and $\phi \in \Omega^n(X)$ in the following way:
\begin{displaymath}
\mu(P)=\int_{N(P)} \omega+\int_P \phi.
\end{displaymath}

The support of a valuation is defined in the obvious way. If $\mu$
is compactly supported, then by \cite{ale05d}, Lemma 2.1.1, one can
choose $\omega$ and $\phi$ to be compactly supported as well. The
integration functional
\begin{displaymath}
\int: \mathcal{V}^\infty_c(X) \to \mathbb{C}
\end{displaymath}
is defined by $\int \mu:=[\phi] \in H_c^n(X)=\mathbb{C}$. Slightly
oversimplifying, $\int \mu=\mu(X)$. Let $s:\p_X \to \p_X, (x,[\xi])
\mapsto (x,[-\xi])$ be the natural involution on $\p_X$. If a smooth
valuation $\mu$ is represented by $(\omega,\phi)$, then $((-1)^n
s^*\omega,(-1)^n\phi)$ represents a smooth valuation $\sigma \mu$.
The involution $\sigma:\mathcal{V}^\infty(X) \to
\mathcal{V}^\infty(X)$ is called {\it Euler-Verdier involution} (see
\cite{ale05b}, Section 3.3).

The pull-back of a smooth valuation $\mu$ under a diffeomorphism $\phi:X' \to X$ is defined by
\begin{displaymath}
\phi^* \mu(P):=\mu(\phi(P)).
\end{displaymath}

The space $\mathcal{V}^\infty(X)$ carries a natural commutative and associative product. Its construction is involved: first a product on smooth valuations on an affine space is constructed in \cite{ale04} and \cite{ale05a} (using the solution of P. McMullen's conjecture \cite{ale01}). Then it is shown in \cite{ale05c} that it can be extended to smooth valuations on an arbitrary manifold $X$ by using local charts. The hard part of this construction is to show that the result is independent of all choices. In Section \ref{section_prod} we will give a much simpler construction of the product which works on the level of differential forms.

One interesting property of the product is the following Poincar\'e duality: the pairing
\begin{align*}
\mathcal{V}^\infty(X) \times \mathcal{V}^\infty_c(X) & \to \mathbb{C},\\
(\mu_1,\mu_2) & \mapsto \int \mu_1 \cdot \mu_2
\end{align*}
is perfect \cite{ale05d} (see also \cite{be07b} for a simpler
proof). Setting
\begin{displaymath}
\mathcal{V}^{-\infty}(X):=(\mathcal{V}_c^\infty(X))^*
\end{displaymath}
one gets an embedding with dense image
\begin{displaymath}
\Xi_\infty:\mathcal{V}^\infty(X) \to \mathcal{V}^{-\infty}(X).
\end{displaymath}
Elements of $\mathcal{V}^{-\infty}(X)$ are called {\it generalized valuations}. There is a canonical embedding with dense image
\begin{align*}
\Xi_\mathcal{P}:\mathcal{P}(X) & \to \mathcal{V}^{-\infty}(X),\\
P & \mapsto [\mu \mapsto \mu(P)].
\end{align*}

\subsection{The Rumin operator}

Let $M$ be a contact manifold of dimension $2n-1$. Recall that this
means that $M$ is given a smooth distribution of codimension $1$
(i.e. a smooth field of hyperplanes in the tangent bundle) which is
completely non-integrable. More explicitly, locally there exists a
1-form $\alpha$ such that the field of hyperplanes is equal to $\ker
\alpha$ with the property that $\alpha \wedge d\alpha^{n-1} \neq 0$.
This form $\alpha$, which is called contact form, is unique up to
multiplication by a non-vanishing smooth function.

A form $\omega \in \Omega^*(M)$ is called {\it vertical} if it vanishes on the contact distribution. Given a local contact form $\alpha$, $\omega$ is vertical if and only if $\omega=\alpha \wedge \phi$ for some $\phi \in \Omega^*(M)$.

Given $\omega \in \Omega^{n-1}(M)$, there exists a {\it unique} vertical form $\omega' \in \Omega^{n-1}$ such that $d(\omega+\omega')$ is vertical (see \cite{rum94}). We define the projection operator $Q:\Omega^{n-1}(M) \to \Omega^{n-1}(M)$ by setting $Q\omega:=\omega+\omega'$. This operator is a first order differential operator containing vertical forms in its kernel. The {\it Rumin operator} is the second order differential operator
\begin{displaymath}
D:=d \circ Q:\Omega^{n-1}(M) \to \Omega^n(M).
\end{displaymath}

The Rumin operator is the main ingredient in some differential complex, called Rumin-de Rham-complex, whose cohomology is isomorphic to the de Rham cohomology \cite{rum94}. We shall not use this isomorphism in the sequel.

If $X$ is a smooth manifold, the cosphere bundle $\p_X:=\p_+(T^*X)$
is defined as the quotient of the cotangent bundle $T^*X$ with the
zero section removed by the natural action of the multiplicative
group $\mathbb{R}_{>0}$ of positive real numbers. This space $\p_X$
carries a natural contact structure. Let $\pi:\p_X \to X$ denote the
projection map. The link between smooth valuations and the Rumin
operator $D$ is given by the following theorem.

\begin{Theorem} \label{thm_bebr} \cite{bebr07}
The smooth valuation represented by $(\omega,\phi) \in \Omega^{n-1}(\p_X) \times \Omega^n(X)$ is trivial if and only if
\begin{align}
D\omega+\pi^*\phi & = 0;\\
\pi_* \omega & = 0.
\end{align}
\end{Theorem}

The first condition is a second order differential equation, while the second condition is a topological condition.

\subsection{Geometric measure theory}

We refer to \cite{fed69} for all notions from geometric measure theory.

In order to fix the notation, we just recall some of the most important definitions. Given an oriented Lipschitz manifold $M$, the current integration over $M$ will be denoted by $[[M]]$, $\int_M$ or just by $M$.

The boundary of a $k$-dimensional current $T$ is the $k-1$-dimensional current $\partial T$ such that $\partial T(\omega)=T(d\omega)$. For $M$ as above, $\partial [[M]]=[[\partial M]]$. A current $T$ with $\partial T=0$ is called a cycle. The support $\spt T$ is defined in the usual way. The space of $k$-dimensional currents on $X$ is denoted by $\mathcal{D}_k(X)$.

If $f:X \to Y$ is a smooth function (or just a Lipschitz function)
such that $f|_{\spt T}$ is proper, then one can define the
push-forward $f_*T$, which is a current on $Y$.

We will need the following simple fact: if $T$ is a $k$-dimensional current of finite mass on a Riemannian manifold $X$,
then there exists a measure $\|T\|$ on $X$ and a simple $k$-vector field $\vec T$ such that $\| \vec T \|=1$ a.e. and

\begin{displaymath}
T(\omega)=\int_X \langle \omega,\vec T\rangle d\|T\|, \quad \omega \in \Omega^k(X).
\end{displaymath}

In particular, rectifiable currents have such representations. In this case $\vec T$
spans the tangent space of $\spt T$ at $\|T\|$-almost all points.

We recall that an integral current is a rectifiable current whose boundary is also rectifiable. The space of integral currents of degree $k$ is denoted by $\mathcal{I}_k(X)$.

We will use at several places the following lemma.
\begin{Lemma} \label{lemma_constancy_theorem}
Let $X$ be a smooth manifold of dimension $n$, $P \in
\mathcal{P}(X)$ and $N(P) \in \mathcal{I}_{n-1}(\p_X)$ its normal
cycle. Assume that the boundary $\partial P$ is connected. If $T \in
\mathcal{I}_{n-1}(\p_X)$ is a cycle with $\spt T \subset \spt N(P)$,
then there exists an integer $m$ with $T=m N(P)$.
\end{Lemma}

\proof
The support of $N(P)$ is some finite union of open manifolds $M_i$ which is homeomorphic to $\partial P$ and hence to an $n-1$-dimensional sphere. If we orient the $M_i$ correctly, then $N(P)=\sum_i [[M_i]]$.

Let $T \in \mathcal{I}_{n-1}(\p_X)$ be a cycle with $\spt T \subset \cup_i M_i$. Applying the constancy theorem from geometric measure theory (\cite{fed69}, 4.1.31) to the restrictions $T|_{M_i}$, we get $T=\sum_i m_i [[M_i]]$ with some integers $m_i$. Since $\partial T=0$, the $m_i$'s are all equal.
\endproof
\subsection{Wave fronts}

We recall the main properties of the wave front of a generalized
function referring for more details to \cite{hormander} or
\cite{gui-st77}, Ch. VI.

Let $X$ be a smooth manifold (always countable at infinity, in
particular paracompact). Let ${\mathcal E}\to X$ be a finite
dimensional vector bundle; for definiteness we assume that
${\mathcal E}$ is a complex bundle though for real bundles the
theory is exactly the same. We denote by $C^\infty(X,{\mathcal E})$
the space of $C^\infty$-sections of ${\mathcal E}$. It is a
Fr\'echet space with topology of uniform convergence on compact
subsets of $X$ of all partial derivatives. We denote by
$C^\infty_c(X,{\mathcal E})$ the space of $C^\infty$-sections of
${\mathcal E}$ with compact support. Naturally
$C^\infty_c(X,{\mathcal E})$ is a locally convex topological vector
space with topology of (strict) countable inductive limit of
Fr\'echet spaces. The inclusion
\begin{displaymath}
C^\infty_c(X,{\mathcal E})\hookrightarrow C^\infty(X,{\mathcal E})
\end{displaymath}
is continuous and has dense image.

Let us denote by $|\omega_X|$ the line bundle over $X$ of complex
densities; thus the fiber of $|\omega_X|$ over a point $x\in X$ is equal
to the one dimensional space of complex valued Lebesgue measures on
the tangent space $T_xX$.

We have a separately continuous bilinear map
\begin{displaymath}
 C^\infty(X,{\mathcal E})\times C^\infty_c(X,{\mathcal E}^*\otimes|\omega_X|)\to \mathbb{C}
\end{displaymath}
given by $(f,g)\mapsto \int_X \langle f,g\rangle$. This map is a non-degenerate
pairing. In other words the induced map
\begin{displaymath}
C^\infty(X,{\mathcal E})\to (C^\infty_c(X,{\mathcal E}^*\otimes|\omega_X|))^*
\end{displaymath}
is continuous, injective and has a dense image when the target
space is equipped with the weak topology.

\begin{Definition}
The space $(C^\infty_c(X,{\mathcal E}^*\otimes|\omega_X|))^*$ is called the space
of {\itshape generalized sections} of ${\mathcal E}$. It is denoted by
$C^{-\infty}(X,{\mathcal E})$.
\end{Definition}
Thus $C^\infty(X,{\mathcal E})\subset C^{-\infty}(X,{\mathcal E})$. Generalized
sections of the trivial line bundle are called {\itshape generalized
functions}. Generalized sections of ${\mathcal E}=|\omega_X|$ are called {\itshape
generalized densities}.

\begin{Remark}
In our previous notation, the space of generalized functions is
$\mathcal{D}_n(X)$ where $n=\dim X$.
\end{Remark}

A main technical tool in the following will be wave front sets. We will not reproduce here their definition, but summarize only the main properties of it relevant for our
applications.

\begin{Proposition} (\cite{hormander} or \cite{gui-st77}, Ch.VI \S
3).\label{P:wave_front_properties} \\
Let $u\in C^{-\infty}(X,{\mathcal E})$.
\begin{itemize}
 \item[(i)] The wave front $\WF(u)$ is a closed
$\R_{>0}$-invariant subset of $T^*X\backslash\underline{0}$, where $\underline{0}$ denotes the zero-section.
\item[(ii)] $\WF(u)=\emptyset$ if and only if $u$ is infinitely smooth.
\item[(iii)] \begin{displaymath} \WF(u\boxtimes v)\subset\left( \WF(u)\times
\WF(v)\right)\cup \left( \WF(u)\times \underline{0}\right)\cup \left(
\underline{0} \times \WF(v)\right).
\end{displaymath}
\end{itemize}
\end{Proposition}

Let $f\colon Y\to X$ be a smooth map, ${\mathcal E}\to X$ be a vector bundle.
Then one has the obvious pull-back map on smooth sections
\begin{displaymath}
f^*\colon C^\infty(X,{\mathcal E})\to C^\infty(Y,f^*{\mathcal E}).
\end{displaymath}
It turns out that $f^*$ can be extended in a natural way to some
generalized sections of ${\mathcal E}$ satisfying appropriate
assumptions. Now we are going to discuss these assumptions referring
for details to \cite{gui-st77}, Ch. VI \S3.

\def\cme{C^{-\infty}(X,\mathcal{E})}
\def\cmela{C^{-\infty}_\Lambda(X,\mathcal{E})}

Let us fix a closed conic subset $\Lambda \subset
T^*X\backslash\underline{0}$. Let us consider the linear subspace
$C^{-\infty}_\Lambda(X,\mathcal{E})\subset
C^{-\infty}(X,\mathcal{E})$ consisting of generalized sections of
$\mathcal{E}$ with the wave front contained in $\Lambda$. The space
$C^{-\infty}_\Lambda(X,\mathcal{E})$ has certain natural locally
convex topology which is stronger than the weak topology induced
from $\cme$ (notice that if $\Lambda=T^*X\backslash\underline{0}$
then both spaces are equal as topological vector spaces).

\begin{Definition}\label{D:transvesality}
Let $\Lambda\subset T^*X\backslash\underline{0}$ be a closed conic
subset. A smooth map $f\colon Y\to X$ is {\itshape transversal to
$\Lambda$} if for any $\xi\in T^*_{f(y)}X \cap \Lambda$ one has
$df^*_y(\xi)\ne 0$.
\end{Definition}

Note that a submersion is transversal to each $\Lambda \subset
T^*X\backslash\underline{0}$.

\begin{Proposition}\label{P:lifting}
If a smooth map $f\colon Y\to X$ is transversal to a closed conic
subset $\Lambda\subset T^*X\backslash\underline{0}$ then there is a
unique linear sequentially continuous map\footnote{A map of
topological spaces is called sequentially continuous if it maps
every  convergent sequence into convergent one. The pull-back map on
generalized functions, as in this theorem, is sequentially
continuous, but not topologically continuous in general.}, also
called pull-back,
\begin{displaymath}
f^*\colon \cmela\to C^{-\infty}_{df^*\Lambda}(Y,f^*\mathcal{E})
\end{displaymath}
where
\begin{displaymath}
df^*(\Lambda):=\{(y,\eta)\in T^*Y|\, \eta\in
df^*_y(\Lambda|_{f(y)})\}
\end{displaymath}
whose restriction to smooth sections of $\mathcal{E}$ is equal to
the pull-back on smooth sections discussed above.
\end{Proposition}

We denote by $T_A^*M$ the conormal bundle of any smooth submanifold
$A\subset M$.

\begin{Remark}\label{restriction-gen-func}
Let $Y$ be a closed submanifold of $X$, and $f\colon Y\to X$ be the
identity imbedding. Then Proposition \ref{P:lifting} says in
particular that $f^*u$ is defined provided
\begin{displaymath}
 \WF(u)\cap T^*_YX=\emptyset.
\end{displaymath}
\end{Remark}

From Proposition \ref{P:wave_front_properties} and Remark
\ref{restriction-gen-func} one can easily deduce the following
result on product of generalized sections (see \cite{gui-st77}, Ch.
VI \S 3, Proposition 3.10). Below we denote by $s: \p_X \to \p_X$
the involution $s((x,[\xi]))=(x,[-\xi])$, and for a subset $Z\subset
\p_X$ we denote by $Z^s$ the image $s(Z)$.

\begin{Proposition} \label{P:tens-prod-wf}
Let ${\mathcal E}_1,{\mathcal E}_2\to X$ be two vector bundles. Let
$\Lambda_1,\Lambda_2$ be closed conic subsets of $T^*X\backslash
\underline{0}$. Let us assume that
\begin{eqnarray}\label{E:tens-prod-wf}
\Lambda_1\cap \Lambda_2^s=\emptyset.
\end{eqnarray}
Let us define a new subset $\Lambda\subset T^*X\backslash
\underline{0}$ such that for any point $x\in X$
\begin{eqnarray}\label{E:wf-of-tens-prod}
\Lambda|_x = (\Lambda_1|_x+\Lambda_2|_x)\cup \Lambda_1|_x\cup
\Lambda_2|_x.
\end{eqnarray}
Then $\Lambda$ is also a closed conic subset, and moreover there
exists a unique bilinear jointly sequentially continuous map
\begin{displaymath}
C^{-\infty}_{\Lambda_1}(X,\mathcal{E}_1)\times
C^{-\infty}_{\Lambda_2}(X,\mathcal{E}_2)\to
C^{-\infty}_{\Lambda}(X,\mathcal{E}_1\otimes \mathcal{E}_2)
\end{displaymath}
whose restriction to smooth sections is the tensor product map.
\end{Proposition}

We will need two further technical propositions.

\begin{Proposition} \label{prop_P:1}
Let $A$ be a smooth submanifold of $M$. Let $\tilde M_A$ denote the
oriented blow up of $M$ along $A$ (see Section \ref{sec_blow_up}
below). Let $f\colon \tilde M_A\to M$ be the natural map. Let $T\in
\mathcal{D}(M)$. Then $f^*T \in \mathcal{D}(\tilde M_A)$ is defined
provided
\begin{displaymath}
 T_A^*M\cap \WF(T)=\emptyset.
\end{displaymath}
\end{Proposition}

\proof
Let $x\in \tilde M_A$. If $x\not\in f^{-1}(A)$ then clearly
$f^*T$ is well defined in a neighborhood of $x$. Thus let us assume
that $x\in f^{-1}(A)$. Let $\xi\in \WF(T)|_{f(x)}$. We have to show
that $(df)^*\xi \in T^*_x \tilde M_A$ does not vanish, or, equivalently, that the
restriction of $\xi$ to $df(T_x\tilde M_A)$ does not vanish.

Since $df(T_x\tilde M_A)$ is a $(\dim A+1)$-dimensional subspace of $T_{f(x)}M$
containing $T_{f(x)}A$, the assumption implies that the restriction of
$\xi$ to $T_{f(x)}A$ does not vanish.
\endproof

\begin{Proposition} \label{prop_P:2}
Let $\tilde M_A \stackrel{f}{\longrightarrow} M \stackrel{g}{\longrightarrow} R$ be smooth maps where
$f$ is the oriented blow up map as in Proposition \ref{prop_P:1}, and $g$
is a submersion. Then $(gf)^*T$ is well defined provided the
following condition is satisfied:
\newline
for any $a\in A$ and any $\zeta \in \WF(T)\cap T^*_{g(a)}R$ the
restriction of $\zeta$ to $dg(T_aA)$ does not vanish.
\end{Proposition}
\proof Since $g$ is a submersion, $g^*T$ is defined, and by
Proposition \ref{P:lifting}
\begin{eqnarray}\label{E:3}
\WF(g^*T)|_a\subset (dg_a)^*(\WF(T)|_{g(a)}).
\end{eqnarray}
By Proposition \ref{prop_P:1} $f^*(g^*T)$ is defined provided
\begin{displaymath}
T^*_AM\cap \WF(g^*T)=\emptyset.
\end{displaymath}

By \eqref{E:3} this condition is satisfied
provided that for any $a\in A$
\begin{displaymath}
(T^*_AM)|_a\cap (dg_a)^*(\WF(T)|_{g(a)})=\emptyset.
\end{displaymath}

The last condition is equivalent to the assumption of the
proposition.
\endproof

\def\dens{|\omega_X|}

Now let us discuss the push-forward of generalized sections. Let
$f\colon Y\to X$ be a smooth {\itshape proper} map. Let
$\mathcal{E}\to X$ be a smooth vector bundle as above. One can
define the push-forward map
\begin{displaymath}
 f_*\colon C^{-\infty}(Y,f^*(\mathcal{E}\otimes \dens^*)\otimes
|\omega_Y|)\to C^{-\infty}(X,\mathcal{E})
\end{displaymath}
as the dual map to
\begin{displaymath}
 f^*\colon C^\infty_c(X,\mathcal{E}^*\otimes\dens)\to
C^\infty_c(Y,f^*(\mathcal{E}^*\otimes\dens)).
\end{displaymath}

Note that $f^*$
indeed takes compactly supported sections to compactly supported
ones due to the properness of our map $f$.

\begin{Remark}\label{R:push-forward-dens}
Let us take $\mathcal{E}=\dens$. Then we get the push-forward map on
generalized densities:
\begin{displaymath}
 f_*\colon C^{-\infty}(Y,|\omega_Y|)\to C^{-\infty}(X,\dens).
\end{displaymath}

In the case when $f$ is a proper submersion, $f_*$ is integration
along the fibers.
\end{Remark}

For a closed conic subset $W\subset T^*Y\backslash\underline{0}$ let
us define a new conic subset
\begin{displaymath}
 df_*W:=\{(x,\eta)\in T^*X \setminus \underline{0}
\,|\,\exists y\in f^{-1}(x) \mbox{ such that } (y,df_y^*(\eta))\in
W|_x\cup\{0\}\}.
\end{displaymath}

Since $f$ is proper, $df_*W$ is closed.

One has the following result (see \cite{gui-st77}, Ch. VI \S 3,
Proposition 3.9).
\begin{Proposition}\label{P:wf-push-forward}
Let $W\subset T^*Y\backslash\underline{0}$ be a closed conic subset.
Then
\begin{displaymath}
f_*\colon C^{-\infty}_W(Y,f^*(\mathcal{E}\otimes
\dens^*)\otimes |\omega_Y|)\to C_{df_*W}^{-\infty}(X,\mathcal{E})
\end{displaymath}
is a sequentially continuous linear map.
\end{Proposition}

Later on we will need the following technical proposition.
\begin{Proposition}\label{P:prop-on-wave-fronts}
Let $X$ be a smooth manifold. Let a Lie group $G$ act smoothly and
transitively on $X$. Let $\mathcal{E}\to X$ be a $G$-equivariant
vector bundle. Let $\Gamma\subset T^*X\backslash\underline{0}_X$ be
a closed conic subset. Let $u\in C^{-\infty}_\Gamma(X,\mathcal{E})$.
Let $\{\mu_j\}$ be a sequence of smooth compactly supported measures
on $G$ with integral 1 and whose supports converge to the identity
element $e\in G$. Then
\begin{displaymath}
u_j:=\int_G(g^*u)\cdot d\mu_j(g)
\end{displaymath}
are smooth sections of $\mathcal{E}$ which converge to $u$ in
$C^{-\infty}_\Gamma(X,\mathcal{E})$ as $j\to \infty$.
\end{Proposition}

\proof The smoothness of $u_j$ follows from the transitivity of the
action of $G$ on $X$. Next, for simplicity of notation we will
assume that $\mathcal{E}$ is the trivial $G$-equivariant line
bundle. Thus $u\in C^{-\infty}_\Gamma(X)$ is a generalized function.
Multiplying $u$ by a smooth compactly supported function, we may
assume that $u$ has a compact support.

Let us fix an arbitrary closed conic neighborhood $\Gamma_1$ of
$\Gamma$ in $T^*X\backslash\underline{0}_X$. Next we can find a
compact symmetric neighborhood $K$ of $e\in G$ such that
$K(\Gamma\cap \pi_X^{-1}(supp(u)))\subset \Gamma_1$, where the
action of $G$ (and hence of $K$) on $T^*X$ is the induced one, and
$\pi_X\colon T^*X\to X$ it the natural map.

Let $p_X\colon G\times X\to X,\, p_G\colon G\times X\to G$ be the
obvious projections. Let $m\colon G\times X\to X$ be the action map.
Let us denote
\begin{displaymath}
\tilde\mu_j:=p_G^*\mu_j\in C^\infty(G\times
X,p^*_G|\omega_G|).
\end{displaymath}

Also we will consider only large $j$ such
that
\begin{displaymath}
 supp(\mu_j)\subset K.
\end{displaymath}

Now let us observe the following identity
\begin{displaymath}
 u_j:=\int_G(g^*u) \cdot d\mu_j(g)=p_{X*}(\tilde \mu_j\cdot m^*u).
\end{displaymath}

It is easy to see that
\begin{displaymath}
 \mu_j \to \delta_e \mbox{ in } C^{-\infty}_{T^*_eG\setminus \underline{0}_G}(G,|\omega_G|).
\end{displaymath}

Hence by Proposition \ref{P:lifting}
\begin{displaymath}
 \tilde\mu_j\to p_G^*\delta_e \mbox{ in }
C^{-\infty}_{(T^*_eG \setminus \underline{0}_G)\times
\underline{0}_X}(G\times X, p_G^*|\omega_G|).
\end{displaymath}

But it is easy to see that $(T^*_eG\times \underline{0}_X)\cap
(m^*\Gamma)=\emptyset$. Then by Proposition \ref{P:tens-prod-wf} we
have a well defined product
\begin{eqnarray}\label{E:product-gener}
C^{-\infty}_{(T^*_eG\backslash\underline{0}_G)\times
\underline{0}_X}(X\times G,p_G^*|\omega_G|)\times
C^{-\infty}_{m^*\Gamma}(X\times G)\to C^{-\infty}_\Lambda(X\times
G,p^*_G|\omega_G|)
\end{eqnarray}
where $\Lambda=\{(y,\xi+\eta)|\,\, (y,\eta)\in m^*\Gamma \mbox{ or }
\eta=0; (y,\xi)\in T^*_eG\times\underline{0}_X\}$. Also the product
(\ref{E:product-gener}) is a jointly sequentially continuous
bilinear map again by Proposition \ref{P:tens-prod-wf}.

Next the push-forward map $p_{X*}\colon C^{-\infty}_\Lambda(X\times
G,p^*_G|\omega_G|)\to C^{-\infty}_{\Gamma_1}(X)$ is sequentially
continuous by Proposition \ref{P:wf-push-forward} since one can
easily see that $dp_{X*}(\Lambda)\subset \Gamma_1$. Thus finally we
deduce that for fixed $u$
\begin{displaymath}
 u_j=p_{X*}(\tilde\mu_j\cdot m^*u)\to u\mbox{ in }
C^{-\infty}_{\Gamma_1}(X).
\end{displaymath}

Since $\Gamma_1$ was an arbitrary
closed conic neighborhood of $\Gamma$, $u_j\to u$ in
$C^{-\infty}_\Gamma(X)$. The proposition is proved.
\endproof

\subsection{Manifolds with corners}\label{Ss:corners}

\begin{Definition} \label{D:corners}
A closed subset $P$ of an $n$-dimensional smooth manifold $X$ is
called a {\itshape submanifold with corners} of dimension $k$ if any
point $p\in P$ has an open neighborhood $U\ni p$ and a
$C^\infty$-diffeomorphism $\phi\colon U\stackrel{\sim}{\to} \R^n$ such that
for some $r\geq 0$
\begin{eqnarray*}
\phi(p)=0,\\
\phi(P\cap U)=\R^{k-r}_{\geq 0}\times \R^r\times 0_{\R^{n-k}}
\end{eqnarray*}
where $0_{\R^{n-k}}$ is the zero element of $\R^{n-k}$. The set of submanifolds with corners is denoted by $\mathcal{P}(X)$.
\end{Definition}

The number $r$ is defined uniquely by the point
$p$, but may depend on it. This $r$ is called {\itshape type} of a
point $p$.

\begin{Example}
\begin{enumerate}
\item Any smooth submanifold of $X$ with or without boundary is a
submanifold with corners.
\item A convex compact $n$-dimensional polytope $P\subset \R^n$ is a
submanifold with corners if and only if $P$ is simplicial, namely
every vertex has exactly $n$ adjacent edges. The vertices are
precisely the points of type $0$.
\end{enumerate}
\end{Example}

A submanifold with corners $P\subset X$ has a natural finite
stratification by locally closed smooth submanifolds as follows. For
any integer $r$, $0\leq r\leq n$ let us denote by $S_r(P)$ the union
of all points of type $r$. Then the $S_r(P)$ are
locally closed smooth disjoint submanifolds of $X$ and
\begin{displaymath}
P=\bigsqcup_{r=0}^{\dim P} S_r(P).
\end{displaymath}

This stratification of $P$ will be called {\itshape canonical
stratification}.

\begin{Definition}\label{D:transversality}
Let $P$ and $Q$ be two closed submanifolds with corners of $X$. We
say that $P$ and $Q$ {\itshape intersect transversally} if each
stratum of the canonical stratification of $P$ is transversal to
each stratum of the canonical stratification of $Q$.
\end{Definition}

\subsection{Double fibrations and Gelfand transform}
\label{section_gelfand_transform}

Let $\pi:M \to A$ be a smooth fiber bundle between oriented manifolds $M$ and $A$ with compact fibers. The fiber integration $(\pi_A)_*:
\Omega^*(M) \rightarrow \Omega^*(A)$ decreases the degree of a form
$\mu$ by the
dimension of the fiber, i.e. by $l:=\dim M-\dim A$. It is defined by
\begin{displaymath}
\int_A \alpha \wedge (\pi_A)_*\mu = \int_M \pi_A^* \alpha \wedge \mu
\end{displaymath}
for all compactly supported differential forms $\alpha$ on $A$. Other sign
conventions can be found in the literature (e.g. \cite{begeve92}); the above one corresponds to the one in \cite{alfe05}. Note that the fiber integration changes its sign if we change the orientation of $M$ or of $A$.

It is easily checked that if $M$ has no boundary then
\begin{equation}
d ((\pi_A)_* \mu)=(\pi_A)_* d\mu
\end{equation}
and that for a form $\alpha$ on $A$ the following projection formula holds:

\begin{equation} \label{projection_formula}
(\pi_A)_* (\pi_A^* \alpha \wedge \mu)=\alpha \wedge (\pi_A)_*
\mu.
\end{equation}

\begin{Lemma} \label{lemma_push_forward_fiber_product}
Let
\begin{displaymath}
\xymatrix{N_1 \ar[d]_{\rho_1}& & N_2 \ar[d]^{\rho_2} \\M_1 \ar[dr]_{\pi_1} & & M_2 \ar[dl]^{\pi_2}\\ & A & }
\end{displaymath}
be a diagram of smooth oriented fiber bundles.
We set
\begin{displaymath}
M_1 \times_A M_2:=\{(m_1,m_2):\pi_1(m_1)=\pi_2(m_2)\}
\end{displaymath}
and
\begin{displaymath}
N_1 \times_A N_2:=\{(n_1,n_2):\pi_1 \circ \rho_1(n_1)=\pi_2 \circ \rho_2(n_2)\}.
\end{displaymath}
Let $p_1,p_2$ be the natural projections from $M_1 \times_A M_2$ to $M_1$ and $M_2$; let $q_1,q_2$ be the natural projections from $N_1 \times_A N_2$ to $N_1$ and $N_2$.
Then for $\nu_1 \in \Omega^{k_1}(N_1)$ and $\nu_2 \in \Omega^{k_2}(N_2)$ we have
\begin{displaymath}
(\rho_1 \times \rho_2)_* q_1^* \nu_1 \wedge q_2^* \nu_2= (-1)^{(\dim N_1+\dim M_1) (k_2+\dim N_2+\dim A)} p_1^* (\rho_1)_*\nu_1 \wedge p_2^* (\rho_2)_* \nu_2.
\end{displaymath}
\end{Lemma}

\begin{Definition}
A double fibration is a diagram
\begin{displaymath}
\xymatrix{& M \ar[dl]_{\pi_A} \ar[dr]^{\pi_B} & \\
A & & B}
\end{displaymath}
where
\begin{enumerate}
\item $\pi_A:M \rightarrow A$ and $\pi_B:M \rightarrow B$ are
  smooth fiber bundles;
\item $\pi_A \times \pi_B:M \rightarrow A \times B$ is a smooth
  embedding.
\end{enumerate}
\end{Definition}

The {\it Gelfand transform} of a differential form $\beta$ on $B$ is the
form $\GT(\beta):=(\pi_A)_* \pi_B^* \beta$.

\begin{Lemma} \label{lemma_push_forwards}
Let
\begin{displaymath}
\xymatrix{M \ar[d]^{\rho_M} \ar[r]^-{\pi_B} & B \ar[d]^{\rho_B}\\
M' \ar[r]^-{\pi_{B'}} & B'}
\end{displaymath}
be a cartesian square (i.e. $M \simeq M' \times_{B'} B$ as oriented manifolds). Then for $\beta \in \Omega^*(B)$
\begin{displaymath}
(\rho_M)_* \circ \pi_B^* \beta=\pi_{B'}^* \circ (\rho_B)_* \beta.
\end{displaymath}
\end{Lemma}

From the lemma, one gets the following functorial property of the Gelfand
transform (\cite{alfe98}, Thm. 2.2).

Let
\begin{displaymath}
\xymatrixcolsep{3pc}
\xymatrix{A \ar[d]^{\rho_A} & M \ar[d]^{\rho_M} \ar[r]^-{\pi_B} \ar[l]_-{\pi_A} & B \ar[d]^{\rho_B}\\
A' & M' \ar[l]_{\pi_{A'}} \ar[r]^-{\pi_{B'}} & B'}
\end{displaymath}
be a morphism of double fibrations (i.e. a commutative diagram of fibrations) such that the right hand part of the diagram is a cartesian square. Then, with $\GT=(\pi_A)_* \circ \pi_B^*$ and $\GT'=(\pi_{A'})_* \circ \pi_{B'}^*$ we have for $\beta \in \Omega^*(B)$
\begin{equation} \label{eq_gelfand_functorial}
(\rho_A)_* \GT(\beta)=\GT'((\rho_B)_* \beta).
\end{equation}

We now suppose that the left hand part of the above diagram is a cartesian square. Then for $\beta' \in \Omega^*(B')$
\begin{equation} \label{eq_gelfand_functorial2}
\GT(\rho_B^* \beta')= \rho_A^* \GT'(\beta').
\end{equation}

Given a current $T$ in $A$, the current $\pi_A^*T$ on $M$ defined by
\begin{displaymath}
\pi_A^*T(\omega):=T((\pi_A)_* \omega)
\end{displaymath}
is called the {\it lift} of $T$ and was studied by Brothers \cite{bro66} and Fu
\cite{fu90}.
In the case of a product bundle $M=A \times F$, the lift
of $T$ is simply $T \times [[F]]$. Moreover, lifting currents is
natural with respect to bundle operations, increases the dimension by
the dimension of the fiber and commutes with the boundary operator $\partial$ if $F$ is without boundary.

\begin{Definition}
The {\it Gelfand transform} of $T$ is
the current $\GT(T):=(\pi_B)_* \pi_A^* T$ in $B$.
\end{Definition}

\section{Blow-up space}
\label{sec_blow_up}

Let $X,Y$ be smooth oriented manifolds of dimensions $n$ and $m$ respectively. We set
\begin{displaymath}
\p_X:=\p_+(T^*X), \p_Y:=\p_+(T^*Y), \p:=\p_+(T^*(X \times Y)).
\end{displaymath}

$\p_X$ consists of pairs $(x,[\xi])$ where $x \in X$, $\xi \in
T_x^*X \setminus \underline{0}$, the brackets mean equivalence class
with respect to the action of $\R^+$ on $T^*X \setminus
\underline{0}$, and $\underline{0}$ means the zero section of the
cotangent bundle. Similarly, $\p_Y$ consists of pairs $(y,[\eta])$
with $y \in Y$ and $\eta \in T_y^*Y \setminus \{0\}$. Finally, $\p$
consists of triples $(x,y,[\xi:\eta])$ with $x \in X, y \in Y$, $\xi
\in T_x^*X, \eta \in T_y^*Y$ not both equal to zero. These manifolds
are canonically oriented.

We define canonical projection maps by the following diagram.
\begin{displaymath}
\xymatrixcolsep{4pc}
\xymatrix{\p_X \ar[d]^{\pi_X} & \p_X \times \p_Y \ar[l]_{q_X} \ar[r]^{q_Y}& \p_Y \ar[d]^{\pi_Y}\\
X & X \times Y \ar[l]_{p_X} \ar[r]^{p_Y} & Y\\
& \p \ar[u]_{\pi_{X \times Y}} &
 }
\end{displaymath}

Let $\mathcal{M}_1 \cong \p_X \times Y \subset \p$ be the
submanifold consisting of triples $(x,y,[\xi:0]) \in \p$. Similarly,
let $\mathcal{M}_2 \cong X\times \p_Y$ consisting of triples
$(x,y,[0:\eta])$.

We set
\begin{displaymath}
\mathcal{M}:=\mathcal{M}_1 \cup \mathcal{M}_2;
\end{displaymath}
clearly the union is disjoint.

Now we are going to describe the oriented blow up of $\mathbb{P}$
along $\mathcal{M}$. Let us consider the fiber bundle
\begin{displaymath}
\p_0:=\p \times_{X \times Y} (\p_X \times \p_Y)
\end{displaymath}
above $X \times Y$. In other words, $\p_0$ is the set of $5$-tuples $(x,y,[\xi:\eta],[\xi'],[\eta'])$, where $x \in X$, $y \in Y$, $[\xi:\eta] \in \p_+(T^*_{(x,y)}(X \times Y))$, $[\xi'] \in \p_+(T_x^*X)$ and $[\eta'] \in \p_+(T^*_yY)$.

Set
\begin{displaymath}
\p_0^1:=\left\{(x,y,[\xi:\eta],[\xi],[\eta]): (x,y,[\xi:\eta]) \in \p \setminus (\mathcal{M}_1 \cup \mathcal{M}_2)\right\} \subset \p_0
\end{displaymath}
and $\hat \p:=\cl(\p_0^1)$, where the closure is taken inside $\p_0$.

The restriction of the natural projection $\p_0 \to \p$ to $\hat \p$ is denoted by
\begin{displaymath}
L: \hat \p \to \p.
\end{displaymath}
Note that $L|_{\p_0^1}: \p_0^1 \to \p \setminus (\mathcal{M}_1 \cup \mathcal{M}_2)$ is a diffeomorphism. We orient $\hat \p$ in such a way that $L|_{\p_0^1}$ is orientation preserving.

We set
\begin{displaymath}
\mathcal{N}_i := L^{-1}(\mathcal{M}_i), \quad i=1,2
\end{displaymath}
and
\begin{displaymath}
\mathcal{N}:=\mathcal{N}_1 \cup \mathcal{N}_2.
\end{displaymath}

Note that
\begin{align*}
\mathcal{N}_1 & = \left\{(x,y,[\xi:0],[\xi],[\eta]): (x,[\xi]) \in \p_X, (y,[\eta]) \in \p_Y\right\};\\
\mathcal{N}_2 & = \left\{(x,y,[0:\eta],[\xi],[\eta]): (x,[\xi]) \in \p_X, (y,[\eta]) \in \p_Y\right\}.
\end{align*}
In particular, there are natural diffeomorphisms $\tau_i:\mathcal{N}_i \to \p_X \times \p_Y, i=1,2$ and
\begin{displaymath}
\dim \mathcal{N}_1=\dim \mathcal{N}_2=2(n+m)-2.
\end{displaymath}

The following diagram commutes.

\begin{displaymath}
\xymatrixcolsep{4pc}
\xymatrix{\mathcal{N}_1 \ar@{^{(}->}[rr] \ar[d]^{\tau_1} & & \hat \p \ar[d]^{L}\\
\p_X \times \p_Y \ar[r]^{id \times \pi_Y}& \p_X \times Y \cong \mathcal{M}_1 \ar@{^{(}->}[r] & \p.
 }
\end{displaymath}

\begin{Lemma}
$\hat \p$ is a $2(n+m)-1$ manifold with boundary $\mathcal{N}$.
\end{Lemma}

\proof
Take a sequence of points $(x_i,y_i,[\xi_i:\eta_i],[\xi_i],[\eta_i])$ in $\p_0^1$, converging to a point
$(x,y,[\xi:\eta],[\xi'],[\eta']) \in \p_0 \setminus \p_0^1$. Then either ($\xi=0$ and $\eta'=\eta$) or ($\eta=0$ and $\xi'=\xi$).
In the first case, the limit point is in $\mathcal{N}_1$, in the second case it is in $\mathcal{N}_2$. Therefore,
\begin{displaymath}
\p_0=\p_0^1 \cup \mathcal{N}.
\end{displaymath}

Let us next show that a neighborhood of $\mathcal{N}_1$ is
diffeomorphic to an open subset of $\p_X \times \p_Y \times \R_{\geq
0}$. For this, we choose Riemannian metrics on $X$ and $Y$. They
induce metrics on the cotangent bundles. A neighborhood of
$\mathcal{N}_1$ in $\p_0$ consists of points of the form
$(x,y,[\xi:\eta'],[\xi],[\eta])$ with $\xi,\eta \neq 0$ and $\eta'
=\lambda \eta$ for some $\lambda \geq 0$. Sending such a point to
$\left(x,[\xi],y,[\eta],\frac{\|\eta'\|}{\|\xi\|}\right) \in \p_X
\times \p_Y \times \R_{\geq 0}$ gives the diffeomorphism we looked
for. The image of $\mathcal{N}_1$ under this diffeomorphism is $\p_X
\times \p_Y \times \{0\}$.

One can prove in a similar way that a neighborhood of
$\mathcal{N}_2$ is diffeomorphic to $\p_X \times \p_Y \times
\R_{\geq 0}$, with $\mathcal{N}_2$ being sent to $\p_X \times \p_Y
\times \{0\}$.
\endproof

Note that $\mathcal{N}=\partial \hat \p$ inherits an orientation from the orientation of $\hat \p$ and that $\tau_i$ is orientation preserving.

We define a map $\Phi$ by
\begin{align*}
\Phi: \hat \p & \to \p_X \times \p_Y,\\
(x,y,[\xi:\eta],[\xi'],[\eta']) & \mapsto ((x,[\xi']),(y,[\eta'])).
\end{align*}

The restriction of $\Phi$ to $\mathcal{N}_i$ is the diffeomorphism $\tau_i$.


\def\cp{\mathcal{P}}
\section{Transversal subsets in $X \times Y$}
\label{sec_transversal} Let $Z$ be a smooth manifold. Let
$P\in\cp(Z)$. Let us have few remarks on the structure of the normal
cycle $N(P)$. First notice that in general $N(P)$ is not a smooth
submanifold (even not with corners), but it is a Lipschitz
submanifold which can be "stratified" in a nice way as explained
below.

As in Section \ref{Ss:corners} we denote by $S_r(P)$ the subset of
$P$ of points of type $r$. We have
$$P=\bigsqcup_{r=0}^{\dim P}S_r(P).$$
Each $S_r(P)$ is a locally closed smooth submanifold of $P$ of
dimension $r$ (without corners or boundary). Let us represent it as
a union of its connected components
$$S_r(P)=S_r^1(P)\bigsqcup\dots\bigsqcup S_r^{l_r}(P).$$
Let us denote the closure
$$N_r^j(P):=\overline{N(P)\cap \pi_Z^{-1}(S_r^j(P))}$$
where $\pi_Z\colon \p_Z\to Z$ is the natural projection as usual.
Clearly
$$N(P)=\bigcup_{r=0}^{\dim P}\bigcup_{j=1}^{l_r}N_r^j(P).$$
Moreover any $N_r^j(P)\subset \p_Z$ is a compact submanifold with
corners of dimension $\dim N_r^j(P)=\dim Z-1(=\dim N(P))$. They have
pairwise disjoint relative interiors.


Let us consider all the strata of the canonical stratification of
all $N_r^j(P)$'s. Let $\{\mathcal{S}_l\}_l$ denote the (finite)
collection of all connected components of all these strata. Thus
$N(P)$ is a disjoint union
$$N(P)=\bigsqcup_l\mathcal{S}_l.$$
Moreover each $N_r^j(P)$ is a disjoint union of some subfamily of
the $\mathcal{S}_l$'s.

Notice that $\mathcal{S}_l$ are locally closed smooth submanifolds
of $\p_Z$ whose closures $\overline{\mathcal{S}_l}$ are compact
submanifolds with corners. Every stratum of the canonical
stratification of $\overline{\mathcal{S}_l}$ is equal to a union of
some other $\mathcal{S}_j$'s.

We can define the smooth part $N^{sm}(P)$ of $N(P)$ by setting
$$N^{sm}(P):=\bigsqcup_{\dim\mathcal{S}_l=\dim N(P)}
\mathcal{S}_l.$$ $N^{sm}(P)$ is open and dense in $N(P)$. The
complement $N(P)\backslash N^{sm}(P)$ has codimension $\geq 1$ in
$N(P)$.

\hfill

In what follows we will apply these notions for a product manifold
$Z=X\times Y$.

\begin{Definition}\label{D:transversal-def}
A set $P\in\cp(X\times Y)$ is called transversal if the smooth
submanifold $\mathcal{M}\subset \p_{X\times Y}$ (defined in Section
\ref{sec_blow_up}) intersects transversally all the strata
$\mathcal{S}_l$ of $N(P)$.
\end{Definition}







Let us denote by $\hat{N}(P)$ the closure of the preimage in
$\hat{\p}$ of $N(P)\backslash \mathcal{M}$. Since $P$ is a
transversal set, $\hat{N}(P)$ can be described as follows:
$$\hat{N}(P)=\bigsqcup_l \hat{\mathcal{S}_l}$$
where $\hat{\mathcal{S}_l}$ is the oriented blow up of
$\mathcal{S}_l$ along $\mathcal{M}\cap \mathcal{S}_l$ which is well
defined because $\mathcal{S}_l$ is a locally closed smooth
submanifold of $\p_{X\times Y}$ intersecting $\mathcal{M}$
transversally; also there is a natural imbedding
$\hat{\mathcal{S}_l}\subset \hat{\p}$. It is easy to see that
$\hat{\mathcal{S}_l}$ are locally closed submanifolds with boundary
in $\hat{\p}$. Their closures $\overline{\hat{\mathcal{S}_l}}$ are
compact submanifolds with corners. Moreover the intersection of any
number of them is again a compact submanifold with corners, and any
connected component of any stratum of the canonical stratification
of such an intersection is equal to a connected component of a
stratum of the canonical stratification of some single
$\overline{\hat{\mathcal{S}_j}}$.


$\hat{N}(P)\subset \hat{\p}$ is a Lipschitz submanifold with
boundary. In order to describe its boundary let us denote
\begin{displaymath}
 N_i(P):=\hat{N}(P)\cap \mathcal{N}_i,\, i=1,2,
\end{displaymath}
these are oriented Lipschitz manifolds of dimension $n+m-2$.
$N_i(P)$ can be stratified as
$$N_i(P)=\bigsqcup_l(\hat{\mathcal{S}}_l\cap \mathcal{N}_i)$$
where $\hat{\mathcal{S}}_l\cap \mathcal{N}_i$ are locally closed
smooth submanifolds. Define the smooth part of $N_i(P)$ by
$$N_i^{sm}(P):=\bigsqcup_{\dim \mathcal{S}_l=\dim N(P)}(\hat{\mathcal{S}}_l\cap
\mathcal{N}_i).$$ $N_i^{sm}(P)$ is open and dense in $N_i(P)$. The
complement $N_i(P)\backslash N_i^{sm}(P)$ has codimension $\geq 1$
in $N_i(P)$.

\hfill

It is easy to see that
\begin{equation} \label{eq_boundary_current}
\partial \hat N(P)=N_1(P)+N_2(P)
\end{equation}
in particular in the sense of currents.

In the case where $X$ and $Y$ are affine spaces, we can define transversal compact convex sets in an analogous way.

\begin{Definition}
A transversal valuation is a functional $\mu$ which is defined on transversal sets $P \in \mathcal{P}(X \times Y)$ and which has the following valuation property: whenever $P_1,P_2,P_1 \cap P_2, P_1 \cup P_2 \in \mathcal{P}(X \times Y)$ are transversal sets, then
\begin{displaymath}
\mu(P_1 \cup P_2)+\mu(P_1 \cap P_2)=\mu(P_1)+\mu(P_2).
\end{displaymath}
The space of transversal valuations is denoted by $\mathcal{V}_{tr}(X \times Y)$.
\end{Definition}

\begin{Lemma} \label{lemma_n1_horizontal}
Let $\omega_1 \in \Omega^*(\p_X)$ be a vertical form and $\omega_2 \in \Omega^*(\p_Y)$ an arbitrary form. Then, for transversal $P \in \mathcal{P}(X \times Y)$,
\begin{displaymath}
\int_{N_1(P)} \tau_1^*(q_X^* \omega_1 \wedge q_Y^* \omega_2)=0.
\end{displaymath}
In particular, this holds if $\omega_1=\pi_X^* \phi_1$ for some $\phi_1 \in \Omega^n(X)$.
\end{Lemma}

\proof If $v$ is a tangent vector to $N_1^{sm}(P)$, then $dL(v) \in
T\p$ is tangent to $N^{sm}(P)$ and thus horizontal. This implies
that $d(q_X \circ \tau_1)v \in T\p_X$ is horizontal and the
assertion follows.
\endproof

\begin{Lemma} \label{lemma_prop_image_projection}
Let $P \in \mathcal{P}(X \times Y)$ be transversal. Then
\begin{align}
((id \times \pi_Y) \circ \tau_1)_* N_1(P) & =0, \label{eq_image_projection_n1}\\
((\pi_X \times id) \circ \tau_2)_* N_2(P) & =0. \label{eq_image_projection_n2}
\end{align}
\end{Lemma}

\proof Let us prove the first equation. The tangent plane to
$N_1^{sm}(P) \subset \mathcal{N}_1$ at a point
$(x,y,[\xi:0],[\xi],[\eta]) \in \mathcal{N}_1$ is generated by lifts
of tangent vectors at $(x,y,[\xi:0])$ of $N^{sm}(P) \cap
\mathcal{M}_1$ and by vectors which are tangent to the fiber
$L^{-1}(x,y,[\xi:0])$. The latter are in the kernel of $(id \times
\pi_Y) \circ \tau_1$.

If $m>1$, then it follows that $((id \times \pi_Y) \circ \tau_1)_* N_1(P)=0$.

In the case $m=1$, we consider the involution
\begin{align*}
s_Y: \mathcal{N}_1 & \to \mathcal{N}_1, \\
(x,y,[\xi:0],[\xi],[\eta]) & \mapsto (x,y,[\xi:0],[\xi],[-\eta]).
\end{align*}
Then $s_Y$ switches the orientation of the fiber of $L$ and
\begin{displaymath}
(s_Y)_* N_1(P)=-N_1(P).
\end{displaymath}
Applying $((id \times \pi_Y) \circ \tau_1)_*$ to this equation yields
\begin{displaymath}
((id \times \pi_Y) \circ \tau_1)_* N_1(P)=-((id \times \pi_Y) \circ \tau_1)_*N_1(P)
\end{displaymath}
and we are done. The proof of the second equation is similar.
\endproof

Let us assume that $X,Y$ are vector spaces and that $A \in
\mathcal{K}^{sm}(X)$ (the space of compact convex subsets of $X$
with smooth boundary and positive curvature) and $B \in
\mathcal{K}^{sm}(Y)$. Let $h_A:X^* \to \R$ and $h_B:Y^* \to \R$
denote the support functions. They are smooth outside the origin and
$1$-homogeneous. In particular the maps $\xi \mapsto d_\xi h_A$ and
$\eta \mapsto d_\eta h_B$ are homogeneous of degree $0$, hence
well-defined on $\p_X$.

We define $H_A,H_B$ by
\begin{align*}
H_A : \p_X & \to \p_X\\
(x,[\xi]) & \mapsto (x+d_\xi h_A,[\xi]),\\
H_B:\p_Y & \to \p_Y\\
(y,[\eta]) & \mapsto (y+d_\eta h_B,[\eta]).
\end{align*}

Let us introduce the following maps:
\begin{align*}
H : \hat \p & \to \p\\
(x,y,[\xi:\eta],[\xi'],[\eta']) & \mapsto (x+d_{\xi'}h_A,y+d_{\eta'}h_B,[\xi:\eta])\\
\tilde H_1: [0,1] \times \mathcal{N}_1 & \to \mathcal{M}_1\\
(t,x,y,[\xi:0],[\xi],[\eta]) & \mapsto (x+d_\xi h_A,y+t d_\eta h_B,[\xi:0])\\
\tilde H_2: [0,1] \times \mathcal{N}_2 & \to \mathcal{M}_2\\
(t,x,y,[0:\eta],[\xi],[\eta]) & \mapsto (x+t d_\xi h_A,y+d_\eta h_B,[0:\eta]).
\end{align*}

Then we have the following commutative diagram:
\begin{equation} \label{eq_diag_h_phi}
\xymatrixcolsep{6pc}
\xymatrix{\hat \p \ar[r]^H \ar[d]_\Phi & \p \ar[d]_{\pi_{X \times Y}}\\
\p_X \times \p_Y \ar[r]^{\pi_X \circ H_A \times \pi_Y \circ H_B} & X \times Y}
\end{equation}

\begin{Proposition}
Let $K \subset X \times Y$ be a transversal compact convex set. Then for $A \in \mathcal{K}^{sm}(X)$ and $B \in \mathcal{K}^{sm}(Y)$ we have
\begin{equation} \label{eq_normal_cycle_sum}
N(K+A \times B)=H_*\hat N(K)-(\tilde H_1)_* ([0,1] \times N_1(K))-(\tilde H_2)_* ([0,1] \times N_2(K)).
\end{equation}
In particular,
\begin{equation} \label{eq_n_and_hat_n}
N(K)=L_* \hat N(K).
\end{equation}
\end{Proposition}

\proof
Let $T$ be the right hand side of \eqref{eq_normal_cycle_sum}. The composition of $\tilde H_i$ with the inclusion $\mathcal{N}_i \hookrightarrow [0,1] \times \mathcal{N}_i, n \mapsto (1,n)$ equals the restriction of $H$ to $\mathcal{N}_i$ ($i=1,2$).

Let $H_i^s$ for $i=1,2$ and $s=0,1$ be the composition of $\tilde H_i$ with the inclusion $\mathcal{N}_i \hookrightarrow [0,1] \times \mathcal{N}_i, n \mapsto (s,n)$. Then $H_1^0$ is the composition of the map $(id \times \pi_Y) \circ \tau_1$ with the map $\p_X \times Y \to \mathcal{M}_1, (x,[\xi],y) \mapsto (x+d_\xi h_A,y,[\xi:0])$. Lemma \ref{lemma_prop_image_projection} implies that $(H_1^0)_* N_1(K)=0$. Similarly, $(H_2^0)_* N_2(K)=0$.

We have $H_i^1=H|_{\mathcal{N}_i}$. By the homotopy formula for currents, $T$ is a cycle.

Next, we will see that $T$ is Legendrian. Since $T$ is integral, it suffices to see that $T$ annihilates vertical forms.

The current $H_* \hat N(K)$ is the image of the restriction of $N(K)$ to $\p \setminus \mathcal{M}$ under the contactomorphism
\begin{align*}
H_{A \times B}=H \circ \left(L|_{\hat \p \setminus \mathcal{N}}\right)^{-1}:\p \setminus \mathcal{M} & \to \p \setminus \mathcal{M}\\
(x,y,[\xi:\eta]) & \mapsto (x+d_\xi h_A,y+d_\eta h_B,[\xi:\eta]).
\end{align*}
Since $N(K)$ annihilates vertical forms, the same holds true for $H_* \hat N(K)$.

Let us consider the second term. We have a commuting diagram
\begin{displaymath}
\xymatrixcolsep{3pc}
\xymatrix{[0,1] \times \mathcal{N}_1 \ar[r]^{\tilde H_1}  \ar[d]_{q_X \circ \tau_1} & \mathcal{M}_1 \ar[d]\\
\p_X \ar[r]^{H_A} & \p_X}
\end{displaymath}
where the vertical map on the right is given by $(x,y,[\xi:0]) \mapsto (x,[\xi])$.

If $v \in T_{(x,y,[\xi:0],[\xi],[\eta])}\mathcal{N}_1$ is a tangent
vector to $N_1^{sm}(K)$, then $(q_X \circ \tau_1)_* v$ is a
horizontal vector in $T_{(x,[\xi])}\p_X$ (see the proof of Lemma
\ref{lemma_n1_horizontal}). Noting that a tangent vector to
$\mathcal{M}_1$ is horizontal if and only if its image in $\p_X$ is
horizontal and that $H_A:\p_X \to \p_X$ is a contactomorphism, it
follows that $(\tilde H_1)_* v$ is horizontal.

Therefore the image under $\tilde H_1$ of almost each tangent plane
of $N_1^{sm}(K)$ is contained in the horizontal distribution. This,
together with the obvious fact that $(\tilde H_1)_*
\left(\frac{\partial}{\partial t}\right)$ is a horizontal vector,
implies that $(\tilde H_1)_* ([0,1] \times N_1(K))$ annihilates
vertical forms. The last term is treated in a similar way.

In order to prove \eqref{eq_normal_cycle_sum}, it remains to show that $T$ and $N(K+A \times B)$ have the same support function. Let us recall some notation from \cite{be07}. Set $V:=X \times Y$. We may identify $\p_+(V^*)$ with the space of oriented hyperplanes in $V$. Let $\mathcal{B}$ be the oriented fiber bundle over $\p_+(V^*)$ whose fiber over a point $[\xi:\eta]$ is given by the oriented line $V / \ker [\xi:\eta]$.

There is a natural map
\begin{align*}
u:\p & \to \mathcal{B}\\
(x,y,[\xi:\eta]) & \mapsto (x,y) / \ker [\xi:\eta].
\end{align*}

Let $\pi_2:\p \to \p_+(V^*)$ be the natural projection. Note that $\pi_2 \circ H_{A \times B}=\pi_2$.  Let $T_{A \times B}:\mathcal{B} \to \mathcal{B}$ be the translation map
\begin{displaymath}
T_{A \times B}((x,y)/\ker [\xi:\eta])=(x+d_\xi h_A,y+d\eta h_B)/\ker [\xi:\eta].
\end{displaymath}
Then the following diagram commutes
\begin{equation} \label{eq_translation_and_u}
\xymatrixcolsep{3pc}
\xymatrix{\p \setminus \mathcal{M} \ar[r]^{H_{A \times B}}  \ar[d]_u & \p \setminus \mathcal{M} \ar[d]_u\\
\mathcal{B} \ar[r]^{T_{A \times B}} & \mathcal{B}}
\end{equation}

The support function of a compact convex set $K \subset V$ can be considered as a Lipschitz function $h_K:\p_+(V^*) \to \mathcal{I}_0(\mathcal{B})$, where $\mathcal{I}_0(\mathcal{B})$ is the space of $0$-dimensional integral currents on $\mathcal{B}$, i.e. finite sums of delta distributions. More precisely, $[\xi:\eta] \in \p_+(V^*)$ is mapped to the delta distribution supported at $\sup \{u(x,y,[\xi:\eta])| (x,y) \in K\} \in \mathcal{B}$. See \cite{be07} and \cite{be06} for details.

The support function of a compactly supported integral Legendrian cycle $T$ on $\p=\p_+(T^*V)$ is the almost everywhere defined function
\begin{align*}
h_T:\p_+(V^*) & \to \mathcal{I}_0(\mathcal{B})\\
[\xi:\eta] & \mapsto u_* \langle T,\pi_2,[\xi:\eta]\rangle.
\end{align*}
Here, the brackets $\langle , \rangle$ mean taking slice in the current-theoretic sense \cite{fed69}.
If $T=N(K)$, then these constructions are compatible in the sense that $h_T=h_K$.

Let, as before, $T$ be the right hand side of \eqref{eq_normal_cycle_sum}. Let $[\xi:\eta] \in \p_+(V^*)$ be a generic point. In particular, $\xi \neq 0, \eta \neq 0$. Then
\begin{align*}
h_T([\xi:\eta]) & = u_* \langle T,\pi_2,[\xi:\eta]\rangle\\
& = u_* \langle (H_{A \times B})_* N(K),\pi_2,[\xi:\eta]\rangle\\
& = u_* (H_{A \times B})_* \langle N(K),\pi_2,[\xi:\eta]\rangle  \quad \text{ (by \cite{fed69}, Thm. 4.3.2 (7))}\\
& = (T_{A \times B})_* h_K([\xi:\eta]) \quad \text{ (by \eqref{eq_translation_and_u})}\\
& = h_{K+A \times B}([\xi:\eta]),
\end{align*}
where the last equation follows from the fact that support functions are additive with respect to Minkowski addition. We thus get that the support functions of $T$ and of $N(K+A \times B)$ agree almost everywhere, which suffices to deduce that $T=N(K+A \times B)$.

Equation \eqref{eq_n_and_hat_n} follows from \eqref{eq_normal_cycle_sum} by taking $A=B=\{0\}$, in which case $H=L$.
\endproof


\section{Exterior product}
\label{section_ext_prod}

The exterior product of smooth valuations on affine spaces $X,Y$ is constructed in \cite{ale05a}. If $\mu_1 \in \mathcal{V}^{\infty}(X)$ is of the form $\mu_1(K)=\phi_1(K+A)$ with $A \in \mathcal{K}^{sm}(X)$ (i.e. $A$ is a smooth convex body with strictly convex boundary) and $\phi_1$ a smooth density on $X$; and if $\mu_2 \in \mathcal{V}^{\infty}(Y)$ is of the form $\nu(K)=\phi_2(K+B)$ with $B \in \mathcal{K}^{sm}(Y)$ and $\phi_2$ a smooth density on $Y$, then
\begin{displaymath}
 \mu \boxtimes \nu(K):=\phi_1 \times \phi_2(K+A \times B), \quad K \in \mathcal{K}(X \times Y).
\end{displaymath}
This definition was extended in \cite{ale05a}, Lemma 4.1.1, to
smooth valuations on affine spaces.

In this section, we extend this exterior product to smooth
valuations on arbitrary smooth manifolds $X,Y$. The construction is
local, so we will assume for simplicity of notation that $X$ and $Y$
are oriented.

\proof[Proof of Theorem \ref{mthm_ext_prod}]
Uniqueness is clear from Properties \eqref{item_invariance} and \eqref{item_affine_case}. Let us prove existence.

Let $\mu_1$ be a smooth valuation on $X$ which is represented by $(\omega_1,\phi_1) \in \Omega^{n-1}(\p_X) \times \Omega^n(X)$. Let furthermore $\mu_2$ be a smooth valuation on $Y$, represented by $(\omega_2,\phi_2) \in \Omega^{m-1}(\p_Y) \times \Omega^m(Y)$.

Let $q_X,q_Y:\p_X \times \p_Y \to \p_X,\p_Y$ be the projection maps. Set $n:=\dim X, m:=\dim Y$. Recall from Section \ref{section_background} that $Q$ is an operator on $\Omega^{n-1}(\p_X)$ (resp. $\Omega^{m-1}(\p_Y)$) and that $D=d \circ Q$.

Let us define forms in $\Omega^{n+m-1}(\hat \p)$ by
\begin{align}
\gamma_0 & :=\Phi^* (q_X^* Q\omega_1 \wedge q_Y^* (D\omega_2+\pi_Y^* \phi_2)),  \label{eq_def_gamma0}\\
\gamma_0' & :=\Phi^* (q_X^* (D\omega_1+\pi_X^* \phi_1) \wedge q_Y^* Q\omega_2), \\
\kappa & := \Phi^*(q_X^* \pi_X^* \phi_1 \wedge q_Y^* Q\omega_2), \label{eq_def_kappa}\\
\kappa' & := \Phi^*(q_X^* Q\omega_1 \wedge q_Y^* \pi_Y^* \phi_2).
\end{align}

Next, we define forms $\gamma_i,i=1,2$ by
\begin{equation} \label{eq_def_gamma_i}
\gamma_i:=\tau_i^* (q_X^* Q\omega_1 \wedge q_Y^* Q\omega_2) \in \Omega^{n+m-2}(\mathcal{N}_i).
\end{equation}

For transversal $P \in \mathcal{P}(X \times Y)$ we define
\begin{equation} \label{eq_ext_prod}
\mu_1 \boxtimes \mu_2(P) =\int_{\hat N(P)} (\gamma_0+(-1)^n \kappa) +(-1)^n \int_{N_1(P)} \gamma_1+\int_P p_X^* \phi_1 \wedge p_Y^* \phi_2.
\end{equation}

We have to show that this definition is independent of the choices of $\omega_i,\phi_i, i=1,2$. Since vertical forms are in the kernel of $Q$ and $D=d \circ Q$, we may suppose that $Q\omega_i=\omega_i$ and thus $D\omega_i=d\omega_i, i=1,2$.

Let us first prove that the exterior product is symmetric, namely
taking into account the orientations, let us show that the following
equation holds:
\begin{equation} \label{eq_symmetry}
\int_{\hat N(P)} (\gamma_0+(-1)^n \kappa)+(-1)^n \int_{N_1(P)} \gamma_1=\int_{\hat N(P)} ((-1)^n \gamma_0'+\kappa')+(-1)^{n+1} \int_{N_2(P)} \gamma_2.
\end{equation}

Since $P \mapsto \hat N(P), P \mapsto N_1(P), P \mapsto N_2(P)$ are valuations,
it suffices to show this equation for all $P$ which are contained in some product $U \times V$ with $U \subset X$ and $V \subset Y$ contractible.
Choose $\rho_1 \in \Omega^{n-1}(U)$ with $d\rho_1=\phi_1$ on $U$ and $\rho_2 \in \Omega^{m-1}(V)$ with $d\rho_2=\phi_2$ on $V$. Set
\begin{displaymath}
\tilde \gamma_0:=\Phi^*( q_X^*(\omega_1+\pi_X^*\rho_1) \wedge q_Y^*
(\omega_2+\pi_Y^*\rho_2)) \in \Omega^{n+m-2}(\hat P \cap
\Phi^{-1}(\pi_X^{-1}U \times \pi_Y^{-1}V)).
\end{displaymath}

Then
\begin{displaymath}
d\tilde \gamma_0 = \gamma_0'+(-1)^{n-1} \gamma_0+\Phi^* (q_X^* (d\omega_1+\phi_1) \wedge q_Y^* \rho_2) + (-1)^{n+1} \Phi^* (q_X^* \rho_1 \wedge q_Y^* (d\omega_2+\phi_2)).
\end{displaymath}
Integrating this equality over $\hat N(P)$ and using \eqref{eq_boundary_current} and Lemma \ref{lemma_prop_image_projection}, one easily proves \eqref{eq_symmetry}.

Next, we prove that the above product is well-defined, i.e. independent of the choices of $(\omega_1,\phi_1)$ and $(\omega_2,\phi_2)$. Suppose that $(\omega_2,\phi_2)$ represents the zero valuation. By Theorem 1 of \cite{bebr07}, we obtain
\begin{equation}
d\omega_2+\pi_Y^* \phi_2=0 \label{eq_closedness}
\end{equation}
and
\begin{equation} \label{eq_cohomological_part}
(\pi_Y)_* \omega_2=0.
\end{equation}

Trivially, it follows that $\gamma_0=0$.

Working locally as above, me may assume that $\phi_2=d\rho_2$ for
some $\rho_2 \in \Omega^{m-1}(V)$. Then $d(\omega_2+\pi_Y^*
\rho_2)=0$ on $\pi_Y^{-1}(V)$ and, using
(\ref{eq_cohomological_part}), there exists $\tau \in
\Omega^{m-2}(\p_Y \cap \pi_Y^{-1}V)$ with $\omega_2+\pi_Y^*
\rho_2=d\tau$ on $\pi_Y^{-1}(V)$.

For transversal $P \in \mathcal{P}(U \times V)$,
\begin{align*}
\int_{N_1(P)} \gamma_1 & = \int_{N_1(P)} \tau_1^*(q_X^* \omega_1 \wedge q_Y^*d\tau) \quad \text{(by Lemma \ref{lemma_prop_image_projection})}\\
& = (-1)^n \int_{N_1(P)} \tau_1^*(q_X^* d\omega_1 \wedge q_Y^* \tau) \quad (\text{since } \partial N_1(P)=0)\\
& = 0 \quad \text{(by Lemma \ref{lemma_n1_horizontal}).}
\end{align*}

Next, we compute that
\begin{align*}
(-1)^n \int_{\hat N(P)} \kappa+ & \int_P p_X^* \phi_1 \wedge p_Y^* \phi_2 = (-1)^n \int_{\hat N(P)} \kappa+(-1)^n \int_{\partial P} p_X^* \phi_1 \wedge p_Y^* \rho_2\\
& = (-1)^n \int_{\hat N(P)} \kappa+(-1)^n \int_{\hat N(P)} \Phi^*(\pi_X^* \phi_1 \wedge \pi_Y^* \rho_2)\\
&  = (-1)^n \int_{\hat N(P)} \Phi^*(\pi_X^* \phi_1 \wedge \pi_Y^* d\tau)\\
& = \int_{N_1(P)} \tau_1^*(q_X^* \pi_X^* \phi_1 \wedge q_Y^* \tau)\\
& \quad +\int_{N_2(P)} \tau_2^*(q_X^* \pi_X^* \phi_1 \wedge q_Y^* \tau)\\
& = 0 \quad \text{(by Lemma \ref{lemma_n1_horizontal} and \eqref{eq_image_projection_n2})}.
\end{align*}

We deduce that the right hand side of \eqref{eq_ext_prod} is independent of the choice of $(\omega_2,\phi_2)$. Taking into account \eqref{eq_symmetry}, similar arguments show that it is also independent of the choice of $(\omega_1,\phi_1)$.

Property \eqref{item_invariance} follows from the construction. Let
us check Property \eqref{item_affine_case}. Suppose that $X$ and $Y$
are affine spaces, $\mu_1$ is given by $\mu_1(K)=\phi_1(K+A)$ with
$A \in \mathcal{K}^{sm}(X)$ and $\phi_1$ a smooth density on $X$;
$\mu_2(K)=\phi_2(K+B)$ with $B \in \mathcal{K}^{sm}(Y)$ and $\phi_2$
a smooth density on $Y$.

We fix $\rho_1 \in \Omega^{n-1}(X)$ with $d\rho_1=\phi_1$. Similarly, let $\rho_2 \in \Omega^{m-1}(Y)$ with $d\rho_2=\phi_2$.

The normal cycle of $K+A$ is given by
\begin{equation}
N(K+A)=H_A(N(K)).
\end{equation}

Using this equation, and the relation $(\pi_X)_* N(K)=\partial K$ for all compact convex $K \subset X$, we obtain
\begin{align*}
\mu_1(K) & =\phi_1(K+A) \\
& = \int_{K+A} \phi_1\\
& = \int_{\partial (K+A)} \rho_1\\
& = \int_{N(K+A)} \pi_X^* \rho_1\\
& = \int_{N(K)} H_A^* \pi_X^* \rho_1.
\end{align*}

It follows that the form $\tilde \omega_1 := H_A^* \pi_X^* \rho_1 \in \Omega^{n-1}(\p_X)$ represents $\mu_1$. Similarly, $\tilde \omega_2:=H_B^* \pi_Y^* \rho_2 \in \Omega^{m-1}(\p_Y)$ represents $\mu_2$. Note that, with $H_A$ and $H_B$ being contactomorphisms, $d\tilde \omega_i$ is vertical.

Define a homotopy
\begin{align*}
\tilde H_A: [0,1] \times \p_X & \to X\\
(t,x,[\xi]) & \mapsto x+t d_\xi h_A
\end{align*}
between $\pi_X \circ H_A$ and $\pi_X$.

Then
\begin{equation} \label{eq_integration_formula}
[[K+A]]=[[K]]+(\tilde H_A)_* ([0,1] \times N(K)).
\end{equation}

Let $r_X:[0,1] \times \p_X \to \p_X$, $r_Y:[0,1] \times \p_Y \to \p_Y$ and $\hat r_i:[0,1] \times \mathcal{N}_i \to \mathcal{N}_i, i=1,2$ denote projections.

Lemma \ref{lemma_push_forward_fiber_product} (with $A=M_1:=\{0\}, N_1:=[0,1], M_2=N_2:=\p_X, \rho_2=id$) implies that
\begin{displaymath}
 [0,1] \times N(K)=(-1)^n r_X^* N(K),
\end{displaymath}
hence \eqref{eq_integration_formula} implies that
\begin{displaymath}
\mu_1(K)=\int_{K+A} \phi_1=(-1)^n \int_{N(K)} (r_X)_* \tilde H_A^* \phi_1+\int_K \phi_1.
\end{displaymath}

Setting $\omega_1:=(-1)^n (r_X)_* \tilde H_A^* \phi_1$, the pair $(\omega_1,\phi_1) \in \Omega^{n-1}(\p_X) \times \Omega^n(X)$ represents the valuation $\mu_1$. Then
\begin{displaymath}
d\omega_1=d\tilde \omega_1-\pi_X^* \phi_1;
\end{displaymath}
which is vertical.

Similarly, with $\omega_2:=(-1)^m (r_Y)_* \tilde H_B^* \phi_2$, the pair $(\omega_2,\phi_2) \in \Omega^{m-1}(\p_Y) \times \Omega^m(Y)$ represents $\mu_2$ and
\begin{equation} \label{eq_relations_eta}
d\omega_2=d\tilde \omega_2-\pi_Y^* \phi_2
\end{equation}
is vertical.

By definition of the exterior product in the affine case (which is also denoted by $\boxtimes$), we obtain for transversal $K$
\begin{align}
\mu_1 \boxtimes \mu_2(K) & = \phi_1 \times \phi_2(K+A \times B) \nonumber \\
& =  \int_{K+A \times B} p_X^* \phi_1 \wedge p_Y^*\phi_2 \nonumber \\
& = \int_{\partial(K+A \times B)} p_X^* \rho_1 \wedge p_Y^* \phi_2 \nonumber \\
& = \int_{N(K+A \times B)} \pi_{X \times Y}^* (p_X^* \rho_1 \wedge p_Y^* \phi_2)  \nonumber
\end{align}

Using \eqref{eq_normal_cycle_sum}  we thus get
\begin{align}
\mu_1 \boxtimes \mu_2(K) & = \int_{\hat N(K)} H^*  \pi_{X \times Y}^* (p_X^* \rho_1 \wedge p_Y^* \phi_2) \nonumber \\
& \quad - \int_{[0,1] \times N_1(K)} \tilde H_1^* \pi_{X \times Y}^* (p_X^* \rho_1 \wedge p_Y^* \phi_2) \nonumber \\
& \quad - \int_{[0,1] \times N_2(K)} \tilde H_2^* \pi_{X \times Y}^* (p_X^* \rho_1 \wedge p_Y^* \phi_2). \label{eq_main_ext_prod}
\end{align}

Lemma \ref{lemma_push_forward_fiber_product} (with $A=M_1=\{0\},N_1=[0,1], N_2=M_2=\mathcal{N}_i, \rho_2=id, k_2=n+m-2$) implies that
\begin{equation} \label{eq_push_forward_N_i}
 [0,1] \times N_i(K) =(-1)^{n+m} \hat r_i^* N_i(K), \quad i=1,2.
\end{equation}

Next, using \eqref{eq_diag_h_phi} we see that
\begin{align}
H^*  \pi_{X \times Y}^* (p_X^* \rho_1 \wedge p_Y^* \phi_2) & =\Phi^*(q_X^* H_A^* \pi_X^* \rho_1 \wedge q_Y^* H_B^* \pi_Y^* \phi_2) \nonumber \\
& = \Phi^*(q_X^* \tilde \omega_1 \wedge q_Y^* d\tilde \omega_2). \label{eq_hat_term}
\end{align}

We apply \eqref{eq_gelfand_functorial2} and Lemma \ref{lemma_push_forward_fiber_product} (with $A=\{0\}$) to the morphism of double fibrations
\begin{displaymath}
\xymatrixcolsep{4pc}
\xymatrix{\mathcal{N}_1 \ar[d]^{\tau_1} & [0,1] \times \mathcal{N}_1 \ar[l]_{\hat r_1} \ar[d] \ar[r]^-{\tilde H_1} & \mathcal{M}_1 \ar[d]^{\pi_{X \times Y}}\\
\p_X \times \p_Y & \p_X \times ([0,1] \times \p_Y) \ar[l]_{id \times r_Y} \ar[r]^-{(\pi_X \circ H_A) \times \tilde H_B } & X \times Y}
\end{displaymath}
where the vertical map in the middle is given by the composition of $id \times \tau_1$ and the natural map $[0,1] \times \p_X \times \p_Y \to \p_X \times [0,1] \times \p_Y$. The left hand square is a cartesian square, but the given orientation of $\mathcal{N}_1$ is opposite to the induced one. We obtain that
\begin{align}
(\hat r_1)_* \tilde H_1^* \pi_{X \times Y}^* (p_X^* \rho_1 \wedge p_Y^* \phi_2)
& = - \tau_1^*( q_X^* (\pi_X \circ H_A)^* \rho_1 \wedge q_Y^* (r_Y)_* \tilde H_B^* \phi_2) \nonumber \\
& = (-1)^{m+1} \tau_1^*( q_X^* \tilde \omega_1 \wedge q_Y^* \omega_2). \label{eq_n1_term}
\end{align}

Similarly, the diagram
\begin{displaymath}
\xymatrixcolsep{5pc}
\xymatrix{
\mathcal{N}_2 \ar[d]_{\tau_2} &  [0,1] \times \mathcal{N}_2 \ar[l]_{\hat r_2} \ar[r]^{\tilde H_2} \ar[d]_{id \times \pi_2} & \mathcal{M}_2 \ar[d]_{\pi_{X \times Y}}\\
\p_X \times \p_Y &  ([0,1] \times \p_X) \times \p_Y \ar[l]_-{r_X \times id} \ar[r]^-{\tilde H_A \times (\pi_Y \circ H_B)} & X \times Y},
\end{displaymath}
commutes. The left hand square is a cartesian square of {\it oriented} manifolds and therefore

\begin{displaymath}
(\hat r_2)_* \tilde H_2^* \pi_{X \times Y}^* (p_X^* \rho_1 \wedge p_Y^* \phi_2)
= (-1)^{m+1} \tau_2^* (q_X^* (r_X)_* \tilde H_A^* \rho_1 \wedge q_Y^* d\tilde \omega_2).
\end{displaymath}

Since $\partial \hat N(K)=N_1(K)+N_2(K)$, we compute that
\begin{align}
\int_{N_2(K)} (\hat r_2)_* \tilde H_2^* \pi_{X \times Y}^* (p_X^* \rho_1 \wedge p_Y^* \phi_2) & = (-1)^{m+1} \int_{\hat N(K)} \Phi^* (q_X^* d (r_X)_* \tilde H_A^* \rho_1 \wedge q_Y^* d\tilde \omega_2) \nonumber \\
& \quad + (-1)^{m}\int_{N_1(K)} \tau_1^* (q_X^* (r_X)_* \tilde H_A^* \rho_1 \wedge q_Y^* d\tilde \omega_2). \label{eq_n2_term}
\end{align}

Next, it is easily computed that
\begin{align} \label{eq_homotopoy_rho1}
d((r_X)_* \tilde H_A^* \rho_1) & =(r_X)_* \tilde H_A^* \phi_1 +(-1)^{n+1} (\pi_X \circ H_A)^* \rho_1+(-1)^n \pi_X^* \rho_1 \nonumber \\
& = (-1)^n \omega_1 + (-1)^{n+1} \tilde \omega_1+(-1)^n \pi_X^* \rho_1.
\end{align}

From the fact that $N_1(K)$ is closed and from \eqref{eq_relations_eta} and \eqref{eq_homotopoy_rho1}, we deduce that
\begin{equation}
\int_{N_1(K)} \tau_1^* (q_X^* (r_X)_* \tilde H_A^* \rho_1 \wedge q_Y^* d\tilde \omega_2)
=\int_{N_1(K)} \tau_1^*(q_X^* (\tilde \omega_1-\omega_1-\pi_X^*\rho_1) \wedge q_Y^* (\omega_2+\pi_Y^*\rho_2)).
\end{equation}

By \eqref{eq_boundary_current} and Lemma \ref{lemma_prop_image_projection}
\begin{align} \label{eq_stokes_nbar}
\int_{N_1(K)} \tau_1^* (q_X^* \pi_X^* \rho_1 \wedge q_Y^* \omega_2) & = \int_{\partial \hat N(K)} \Phi^*(q_X^* \pi_X^* \rho_1 \wedge q_Y^* \omega_2) \nonumber \\
& = \int_{\hat N(K)} \Phi^* (q_X^* \pi_X^* \phi_1 \wedge q_Y^* \omega_2) \nonumber \\
& \quad +(-1)^{n-1} \int_{\hat N(K)} \Phi^* (q_X^* \pi_X^* \rho_1 \wedge q_Y^* d\omega_2).
\end{align}

We obtain
\begin{align*}
\mu_1 \boxtimes \mu_2(K) & = \int_{\hat N(K)} \Phi^*\left(q_X^* \tilde \omega_1 \wedge q_Y^* d\tilde \omega_2\right)\\
& \quad +(-1)^{n+m+1} \int_{N_1(K)} (r_X)_*\tilde H_1^*\pi_{X \times Y}^*(p_X^*\rho_1 \wedge p_Y^*\phi_2)\\
& \quad +(-1)^{n+m+1} \int_{N_2(K)} (r_Y)_*\tilde H_2^*\pi_{X \times Y}^*(p_X^*\rho_1 \wedge p_Y^*\phi_2)\\
& \quad \quad \text{ (by \eqref{eq_main_ext_prod}, \eqref{eq_hat_term},\eqref{eq_push_forward_N_i})}\\
& = \int_{\hat N(K)} \Phi^*\left(q_X^* \tilde \omega_1 \wedge q_Y^* d\tilde \omega_2\right)\\
& \quad +(-1)^n \int_{N_1(K)} \tau_1^*(q_X^* \tilde \omega_1 \wedge q_Y^* \omega_2)\\
& \quad +(-1)^n  \int_{\hat N(K)} \Phi^*(q_X^*d(r_X)_*\tilde H_A^*\rho_1 \wedge q_Y^* d\tilde \omega_2)\\
& \quad +(-1)^{n+1} \int_{N_1(K)} \tau_1^*(q_X^*(r_X)_*\tilde H_A^*\rho_1 \wedge q_Y^*d\tilde \omega_2) \\
& \quad \quad \text{ (by \eqref{eq_n1_term}, \eqref{eq_n2_term})}.
\end{align*}

Replacing \eqref{eq_relations_eta} and \eqref{eq_homotopoy_rho1} into this equation gives us
\begin{align*}
\mu_1 \boxtimes \mu_2(K) & = \int_{\hat N(K)} \Phi^*(q_X^* \tilde \omega_1 \wedge q_Y^* (d\omega_2+\pi_Y^*\phi_2))\\
& \quad -\int_{\hat N(K)} \Phi^*(q_X^*(-\omega_1 + \tilde \omega_1- \pi_X^* \rho_1) \wedge q_Y^* (d\omega_2+\pi_Y^*\phi_2))\\
& \quad +(-1)^n \int_{N_1(K)} \tau_1^*(q_X^* \tilde \omega_1 \wedge q_Y^* \omega_2)\\
& \quad +(-1)^{n+1} \int_{N_1(K)} \tau_1^*(q_X^* (\tilde \omega_1-\omega_1-\pi_X^*\rho_1) \wedge q_Y^* (\omega_2+\pi_Y^*\rho_2)) \\
\end{align*}

The sum of the first two terms is
\begin{multline*}
 \int_{\hat N(K)} \Phi^*(q_X^*(\omega_1 + \pi_X^* \rho_1) \wedge q_Y^* (d\omega_2+\pi_Y^*\phi_2))\\=\int_{\hat N(K)} \gamma_0+\int_{\hat N(K)} \Phi^*( q_X^*\pi_X^* \rho_1 \wedge q_Y^* d\omega_2)+\int_{\hat N(K)} \Phi^*(q_X^*\pi_X^* \rho_1 \wedge q_Y^* \pi_Y^*\phi_2),
\end{multline*}
while the sum of the last two terms is, by Lemma \ref{lemma_prop_image_projection},
\begin{multline*}
 (-1)^{n} \int_{N_1(K)} \tau_1^*(q_X^*(\omega_1+\pi_X^*\rho_1) \wedge q_Y^* \omega_2)\\
=(-1)^{n} \int_{N_1(K)} \gamma_1 + (-1)^{n} \int_{N_1(K)} \tau_1^*(q_X^* \pi_X^*\rho_1 \wedge q_Y^* \omega_2)
\end{multline*}

From \eqref{eq_image_down_normal_cycle} and \eqref{eq_n_and_hat_n} we infer that
\begin{displaymath}
 (\pi_{X \times Y} \circ L)_*\hat N(K)=(\pi_{X \times Y})_* N(K)=\partial K,
\end{displaymath}
and therefore

\begin{align*}
\int_{\hat N(K)} \Phi^*(q_X^* \pi_X^* \rho_1 \wedge q_Y^* \pi_Y^* \phi_2) & = \int_{\hat N(K)} (\pi_{X \times Y} \circ L)^* (p_X^* \rho_1 \wedge p_Y^* \phi_2)\\
& = \int_{\partial K} p_X^* \rho_1 \wedge p_Y^* \phi_2\\
& = \int_K p_X^* \phi_1 \wedge p_Y^* \phi_2.
\end{align*}

Next, by Lemma \ref{lemma_prop_image_projection} and Stokes's theorem, we have
\begin{align*}
\int_{N_1(K)} \tau_1^*(q_X^*\pi_X^*\rho_1 \wedge q_Y^*\omega_2) & = \int_{\partial \hat N(K)} \Phi^*(q_X^*\pi_X^*\rho_1 \wedge q_Y^*\omega_2)\\
& = \int_{\hat N(K)}d\Phi^*(q_X^*\pi_X^*\rho_1 \wedge q_Y^*\omega_2)\\
& = \int_{\hat N(K)} \kappa +(-1)^ {n-1} \int_{\hat N(K)} \Phi^*( q_X^*\pi_X^* \rho_1 \wedge q_Y^* d\omega_2)
\end{align*}

This gives us
\begin{align*}
 \mu_1 \boxtimes \mu_2(K) & = \int_{\hat N(K)}(\gamma_0+(-1)^n\kappa)+(-1)^n\int_{N_1(K)}\gamma_1+\int_K p_X^*\phi_1 \wedge p_Y^*\phi_2,
\end{align*}
which finishes the proof of Theorem \ref{mthm_ext_prod}.
\endproof


\section{Product}
\label{section_prod}

The product of two smooth valuations is defined to be the restriction of the exterior product to the diagonal in $X \times X$. In this section, we give a description of the product in terms of differential forms.

The following commutative diagram may help to understand the definition of the various maps which will be introduced in this section.
\begin{displaymath}
\xymatrixcolsep{4pc}
\xymatrix{\overline{\mathcal{N}}_i \ar[d]^{\bar p \circ \bar \iota_i} \ar[rd]^{\bar\iota_i}
 & \overline{\mathcal{N}}_0  \ar[d]^{\bar \iota_0} \ar[r]^{(\bar L \circ \bar \iota_0)|_{\bar{\mathcal{N}}_0}} & \mathcal{N}_0 \ar[d]  & \mathcal{N}_i \ar[dl]\\
\p_X \ar[d]^{\pi} & \bar \p \ar[d]^{q_1 \circ \bar \Phi} \ar[r]^{\bar L} \ar[rd]^{\bar \Phi} \ar[l]_{\bar p} & \hat \p \ar[r]^L \ar[d]^{\Phi} & \p \cap \pi_{X \times X}^{-1}\Delta \ar[d]^{p_1 \circ \pi_{X \times X}}\\
X & \p_X \ar[l]_{\pi} \ar[uur]_\rho  & \p_X \times_X \p_X \ar[r]^{\pi \circ q_2}  \ar[l]_{q_1}& X\\
 }
\end{displaymath}

Let $X$ be a smooth oriented $n$-dimensional manifold.

Let us fix the following notation: let $q_1,q_2:\p_X \times_X \p_X \to \p_X$, $\pi:\p_X \to X$ and $\pi_{X \times X}:\p \to X \times X$ be the natural projections. By $\Delta$ we will denote the diagonal embedding from $X$ to $X \times X$ and also the diagonal in $X \times X$. The projections from $X \times X$ to its factors are denoted by $p_1,p_2$.

Let us consider the fiber bundle over $X$ consisting of tuples $(x,[\xi:\eta],[\xi'],[\eta'],[\zeta])$ with $x \in X$, $[\xi:\eta] \in \p_+(T_{(x,x)}^* X \times X)$, $[\xi'] \in \p_+(T_x^*X)$, $[\eta'] \in \p_+(T_x^*X)$ and $[\zeta] \in \p_+(T_x^*X)$.

Let $\bar \p$ be the closure of the set of all such tuples with $\xi \neq 0, \eta \neq 0, \xi+\eta \neq 0$ and $[\xi']=[\xi],[\eta']=[\eta], [\xi+\eta]=[\zeta]$.
Then $\bar \p$ is a compact $3n-1$-dimensional manifold with boundary. More precisely, the boundary consists of three disjoint submanifolds:
\begin{displaymath}
\partial \bar \p= \bar{\mathcal{N}}_0 \cup \bar{\mathcal{N}}_1 \cup \bar{\mathcal{N}}_2 =:\bar{\mathcal{N}}.
\end{displaymath}
Here  $\bar{\mathcal{N}}_0$ is the set of all tuples of the form $(x,[\xi:-\xi],[\xi],[-\xi],[\zeta])$; $\bar{\mathcal{N}_1}$ is the set of all tuples of the form $(x,[\xi:0],[\xi],[\eta'],[\xi])$ and $\bar{\mathcal{N}}_2$ is the set of all tuples of the form $(x,[0:\eta],[\xi'],[\eta],[\eta])$.

We set $\bar{\iota}_i:\bar{\mathcal{N}}_i \hookrightarrow \bar \p, i=0,1,2$ for the embeddings.

Recall the definition of $\hat \p$ from Section \ref{sec_blow_up}. Here we modify the definition by letting $\hat \p$ the closure (inside $\hat \p_0$) of the set
\begin{displaymath}
 \left\{(x,x,[\xi:\eta],[\xi],[\eta])| (x,x,[\xi:\eta]) \in \p \setminus (\mathcal{M}_1 \cup \mathcal{M}_2)\right\}.
\end{displaymath}
In other words, we only take the part of the previously defined $\hat \p$ where $x=y$. Then $\hat \p$ is a $3n-1$-dimensional manifold with boundary. Instead of writing $(x,x,[\xi:\eta],[\xi'],[\eta'])$, we will use the simpler notation $(x,[\xi:\eta],[\xi'],[\eta'])$.

Let
\begin{align*}
\bar L:\bar \p & \to \hat \p,\\
(x,[\xi:\eta],[\xi'],[\eta'],[\zeta]) & \mapsto (x,[\xi:\eta],[\xi'],[\eta'])
\end{align*}
be the natural projection and set $\bar \Phi:=\Phi \circ \bar L:\bar \p \to \p_X \times_X \p_X$.

Note that $\bar{\mathcal{N}}_i=(\bar L \circ \bar \iota_i)^{-1}(\mathcal{N}_i)$ where $\mathcal{N}_1,\mathcal{N}_2$ are defined as before and $\mathcal{N}_0$ is the set of tuples of the form $(x,[\xi:-\xi],[\xi],[-\xi])$.

Define
\begin{align*}
\bar p:\bar \p & \to \p_X,\\
(x,[\xi:\eta],[\xi'],[\eta'],[\zeta]) & \mapsto (x,[\zeta]).
\end{align*}
Note that $\bar p \circ \bar \iota_i:\bar{\mathcal{N}}_i \to \p_X$ is a fiber bundle whose fiber is diffeomorphic to $S^{n-1}$.

There is a commuting diagram of fiber bundles
\begin{equation} \label{eq_orientation_bar_p}
\xymatrix{
\bar \p \ar[d]_-{\pi \circ \bar p} \ar[r]^-{L \circ \bar L} & \p\cap \pi_{X \times X}^{-1}(\Delta)  \ar[d]^-{p_1 \circ \pi_{X \times X}} \\
X \ar[r]^{id} & X}
\end{equation}

We orient $\bar \p$ in such a way that the restriction of $L \circ \bar L$ to each fiber of $\pi \circ \bar p$ is an orientation preserving diffeomorphism to its image (which is a dense open subset of the corresponding fiber of $p_1 \circ \pi_{X \times X}$).

Consider the double fibration
\begin{displaymath}
\xymatrix{
& \bar \p \ar[dl]_{\bar p} \ar[dr]^{\bar \Phi} &  \\
\p_X & & \p_X \times_X \p_X}
\end{displaymath}

and the corresponding Gelfand transform
\begin{align*}
\GT: \Omega^*(\p_X \times_X \p_X) & \to \Omega^{*-n}(\p_X)\\
\omega & \mapsto \bar p_* \bar \Phi^* \omega.
\end{align*}

\begin{Lemma} \label{lemma_boundary_terms_push_forward}
Let $\omega \in \Omega^k(\bar \p)$. Then
\begin{equation}
d\bar p_* \omega=\bar p_* d\omega+(-1)^{n-k-1} \sum_{i=0}^2 (\bar p \circ \bar \iota_i)_* \bar \iota_i^* \omega.
\end{equation}
\end{Lemma}

\proof
Recall that the push-forward (or fiber integration) $\bar p_*:\Omega^{*}(\bar \p) \to \Omega^{*-n}(\p_X)$ is defined by the relation
\begin{displaymath}
\int_{\p_X} \beta \wedge \bar p_* \omega = \int_{\bar \p} \bar p^* \beta \wedge \omega, \quad \forall \beta \in \Omega^*(\p_X).
\end{displaymath}
From this, we obtain for $\beta \in \Omega^{3n-k-2}(\p_X)$
\begin{align*}
\int_{\bar \p} d(\bar p^* \beta \wedge \omega) & =\int_{\bar \p} \bar p^* d\beta \wedge \omega+(-1)^{n-k} \int_{\bar \p} \bar p^* \beta \wedge d\omega\\
& = \int_{\p_X} d\beta \wedge \bar p_*\omega +(-1)^{n-k} \int_{\p_X} \beta \wedge \bar p_* d\omega\\
& = (-1)^{n-k-1} \int_{\p_X} \beta \wedge d \bar p_*\omega +(-1)^{n-k} \int_{\p_X} \beta \wedge \bar p_* d\omega.
\end{align*}

On the other hand, by Stokes's theorem,
\begin{displaymath}
\int_{\bar \p} d(\bar p^* \beta \wedge \omega) = \sum_{i=0}^2 \int_{\bar{\mathcal{N}}_i} \bar \iota_i^*(\bar p^* \beta \wedge \omega) = \sum_{i=0}^2 \int_{\p_X} \beta \wedge  (\bar p \circ \bar \iota_i)_* \bar \iota_i^* \omega.
\end{displaymath}
\endproof

\begin{Lemma} \label{lemma_verticality}
If $\omega_i \in \Omega^*(\p_X), i=1,2$ are vertical forms, then
$\GT(q_1^* \omega_1 \wedge q_2^* \omega_2)$ is vertical.
\end{Lemma}

\proof
First note that
\begin{equation} \label{eq_projections_x}
\pi \circ q_1 \circ \bar \Phi= \pi \circ q_2 \circ \bar \Phi= \pi \circ \bar p.
\end{equation}

Fix a point $(x,[\xi:\eta],[\xi],[\eta],[\zeta]) \in \bar \p \setminus \bar{\mathcal{N}}$ and representatives $\xi,\eta,\zeta$ of the corresponding equivalence classes. Then $\xi+\eta=c \zeta$ for some $c>0$.
Set $\alpha_1:=\pi^* \xi$, $\alpha_2:=\pi^* \eta$, $\alpha:=\pi^* \zeta$. Since $\omega_1$ is vertical, it is (at the point $(x,[\xi])$) divisible by $\alpha_1$. Similarly, $\omega_2$ is divisible by $\alpha_2$.
Applying \eqref{eq_projections_x}, we obtain
\begin{displaymath}
\bar \Phi^* (q_1^* \alpha_1 \wedge q_2^* \alpha_2)=\bar p^* \alpha_1 \wedge c \bar p^* \alpha.
\end{displaymath}

The same argument applies to each point of the fiber $\bar p$ and shows that $\bar \Phi^* (q_1^* \omega_1 \wedge q_2^* \omega_2)$ is divisible by $\bar p^* \alpha$. The push-forward of such a form under $\bar p$ is a multiple of $\alpha$ by the projection formula \eqref{projection_formula}.
\endproof

Let us introduce a map $\bar r_0:=q_1 \circ \bar \Phi \circ \bar \iota_0$ and the diffeomorphisms
\begin{align*}
\rho:\p_X & \to \mathcal{N}_0, \\
(x,[\xi]) & \mapsto (x,[\xi:-\xi],[\xi],[-\xi]);\\
s: \p_X & \to \p_X,\\
(x,[\xi]) & \mapsto (x,[-\xi]).
\end{align*}

 Then we have a commutative diagram
\begin{equation} \label{eq_diag_r0}
\xymatrixcolsep{3pc}
\xymatrix{
\bar{\mathcal{N}}_0 \ar[r]^-{\bar L \circ \bar \iota_0} \ar[d]_{\bar r_0} &\mathcal{N}_0 \ar[d]^{\Phi|_{\mathcal{N}_0}}\\
\p_X \ar[r]_-{(id,s)} \ar[ur]^\rho         & \p_X \times \p_X}
\end{equation}

\begin{Lemma} \label{lemma_normal_cycle_diagonal}
Let $P \in \mathcal{P}(X)$. Then $\Delta(P)\in \mathcal{P}(X\times
X)$ is transversal  in the sense of Definition
\ref{D:transversal-def}. If moreover $P$ is diffeomorphic to a
simplex then
\begin{align}
\hat N(\Delta P) & =\bar L_* \bar p^* N(P)+\rho_* \pi^* P\\
N_i(\Delta P) & = (-1)^{n} (\bar L \circ \bar \iota_i)_* (\bar p \circ \bar \iota_i)^* N(P), i=1,2.
\end{align}
\end{Lemma}

\proof

First let us show the claim about transversality. Since the claim is
local we may and will assume that $X$ is the germ of $\R^n$ at the
point $0$ and
$$P=\R^{k-l}_{\geq 0}\times \R^l\times 0_{n-k}.$$
Then we have to show that $\mathcal{M}$ intersects transversally any
stratum of $S_l(\Delta(P))$. Let $T$ be an arbitrary such stratum.
Up to a linear change of coordinates in $\R^n$, $T$ has the
following form. There exists $m\leq k-l$ such that $T$ consists of
the elements of the form $(x,x,\zeta)$ where $x=(0_{k-l},y,0_{n-k})$
with $y\in \R^l$ is arbitrary, and $\zeta\in \R^n\times \R^n$
belongs to the submanifold $T'$ of the unit sphere of
$\R^n\times\R^n=\R^{2n}$ (here we identify $\R^{n*}$ with $\R^n$ via
the standard pairing) defined by the following two conditions:

(a) $|\zeta|=1$;

$\mbox{(b) }\zeta=(\kappa,\kappa)+(\omega,-\omega)$ where $\omega\in
\R^n$ is arbitrary, $\kappa=(\kappa_m,0_{k-l-m},0_l,\lambda)$ with
$\kappa_m\in \R_{<0}^m,\lambda\in \R^{n-k}$.

By the symmetry it suffices to show that $\mathcal{M}_1$ intersects
$T$ transversally. Since $\mathcal{M}_1$ is invariant under the
translations by $X\times X$, it suffices to show that the manifold
$T'$ intersects transversally in $S^{2n-1}$ with another submanifold
$$M':=\{(\xi,0_n)|\, \xi\in \R^n, |\xi|=1\}.$$


Let us denote by $T''$ and $M''$ the cones over $T'$ and $M'$
respectively in $\R^n\times \R^n$. It suffices (in fact equivalent)
to show that $T''$ and $M''$ intersect transversally everywhere,
except possibly 0. We have
$$M''=\{(\xi,0_n)| \xi\in \R^n\}.$$
On the other hand $T''$ consists of vectors satisfying the condition
(b) above. Clearly $M''$ is a linear subspace of $\R^{2n}$, and
$T''$ is an open convex cone in the linear subspace $\tilde T$
consisting, by definition, of the vectors
$$\zeta=(\kappa,\kappa)+(\omega,-\omega)$$ where $\omega\in \R^n$ is
arbitrary, $\kappa=(\kappa_m,0_{k-l-m},0_l,\lambda)$ with
$\kappa_m\in \R^m,\lambda\in \R^{n-k}$. It suffices to show that the
linear subspaces $M''$ and $\tilde T$ intersect transversally, i.e.
$M''+\tilde T=\R^{2n}$. But this is obvious.

\hfill

Let us prove now the second part of the lemma. Set
\begin{align*}
T & :=\bar L_* \bar p^* N(P)\\
Q & := \rho_* \pi^* P\\
T_i &:=(\bar L \circ \bar \iota_i)_* (\bar p \circ \bar \iota_i)^* N(P), \quad i=0,1,2.
\end{align*}

Clearly $\spt T_i \subset \mathcal{N}_i$ for $i=0,1,2$.

Since $N(P)$ is closed and using Lemma \ref{lemma_boundary_terms_push_forward} we compute that
\begin{equation}
\partial T=(-1)^n (T_0+T_1+T_2). \label{eq_boundary_of_t}
\end{equation}

Tangent vectors to the fibers of $\bar p \circ \bar \iota_i, i=1,2$ are in the kernel of $d(L \circ \bar L \circ \bar \iota_i)$, and hence $L_* T_1=0, L_*T_2=0$. Since $\partial P=\pi_* N(P)$, applying \eqref{eq_gelfand_functorial2} to the morphism of double fibrations
\begin{displaymath}
\xymatrixcolsep{3pc}
\xymatrix{\p_X \ar[d]_{\pi} & \bar{\mathcal{N}}_0 \ar[d]_{\bar r_0} \ar[l]_-{\bar p \circ \bar \iota_0} \ar[r]^-{\bar L \circ \bar \iota_0} & \mathcal{N}_0 \ar[d]_{id}\\
X & \p_X \ar[r]^{\rho} \ar[l]_-{\pi} & \mathcal{N}_0}
\end{displaymath}
yields $T_0=(-1)^{n+1} \partial Q$. It follows that
\begin{displaymath}
\partial L_*(T+Q)=0.
\end{displaymath}
The integral current $L_*(T+Q)$ is given by integration over a submanifold of $\p$. This submanifold is contained in the support of $N(\Delta P)$. Lemma \ref{lemma_constancy_theorem} implies that $L_*(T+Q)=m N(\Delta P)$ for some $m \in \mathbb{Z}$.

In order to fix $m$, we let $P$ shrink to a point $x \in X$. Then $Q$ weakly converges to $0$; $T$ weakly converges to $\bar L_* \bar p^* N(\{x\})$ and $N(\Delta P)$ weakly converges to $N(\{(x,x)\})$. Applying Lemma \ref{lemma_push_forwards} to \eqref{eq_orientation_bar_p} yields
\begin{align*}
L_* \bar L_* \bar p^* N(\{x\}) & = (L \circ \bar L)_* \bar p^* \pi^* \delta_x\\
& = (p_1 \circ \pi_{X \times X})^* \delta_x\\
& = N(\{(x,x)\}).
\end{align*}
We deduce that $m=1$ and thus $\hat N(\Delta P)=T+Q$. Now \eqref{eq_boundary_of_t} implies that $N_i(\Delta P)=(-1)^{n} T_i,i=1,2$.
\endproof

\begin{Proposition} \label{prop_rumin_omega}
Given $(\omega_i,\phi_i) \in \Omega^{n-1}(\p_X) \times \Omega^n(X), i=1,2$, define forms $\omega,\phi \in \Omega^{n-1}(\p_X) \times \Omega^n(X)$ by
\begin{align}
\omega & = \GT(q_1^* \omega_1 \wedge q_2^* D\omega_2) + \omega_1 \wedge \pi^*\pi_* \omega_2 \in \Omega^{n-1}(\p_X), \nonumber\\
\phi & = \pi_* (\omega_1 \wedge s^* (D\omega_2+\pi^* \phi_2))+ \phi_1 \wedge \pi_* \omega_2 \in \Omega^n(X). \label{eq_forms_prod}
\end{align}
Then
\begin{align} \label{eq_Rumin_omega}
D\omega+\pi^*\phi & =\GT(q_1^*(D\omega_1+\pi^*\phi_1) \wedge q_2^* (D\omega_2+\pi^*\phi_2)) \nonumber\\
& \quad +\pi^*\pi_* \omega_1 \wedge (D\omega_2+\pi^*\phi_2) + \pi^*\pi_*\omega_2 \wedge (D\omega_1+\pi^*\phi_1)\\
\pi_*\omega & = \pi_*\omega_1 \wedge \pi_* \omega_2. \label{eq_push_forward_omega}
\end{align}
\end{Proposition}

\proof
By Lemma \ref{lemma_verticality} and the fact that vertical forms are in the kernel of $D$, we may assume that $d\omega_i=D\omega_i,i=1,2$.

We apply \eqref{eq_gelfand_functorial2} to the morphism of double fibrations
\begin{displaymath}
\xymatrixcolsep{3pc}
\xymatrix{\p_X \ar[d]_{\pi} & \bar{\mathcal{N}}_0  \ar[d]^-{\bar r_0} \ar[l]_{\bar p \circ \bar \iota_0} \ar[r]^-{\bar \Phi \circ \bar \iota_0} &  \p_X \times_X \p_X \ar[d]^{id}\\
X & \p_X \ar[l]_\pi \ar[r]^-{(id,s)} & \p_X \times_X \p_X}.
\end{displaymath}
It follows that
\begin{align*}
(\bar p \circ \bar \iota_0)_* \bar \iota_0^* \bar \Phi^* (q_1^* \omega_1 \wedge q_2^* D\omega_2) =(-1)^{n+1} \pi^* \pi_*(\omega_1 \wedge s^*D\omega_2).
\end{align*}

We have a morphism of double fibrations
\begin{equation} \label{eq_double_n1}
\xymatrixcolsep{3pc}
\xymatrix{\p_X \ar[d]^{id} & \bar{\mathcal{N}}_1 \ar[d]^{\bar p \circ \bar \iota_1} \ar[r]^-{\bar \Phi \circ \bar \iota_1} \ar[l]_-{\bar p \circ \bar \iota_1} & \p_X \times_X \p_X \ar[d]^{id \times \pi}\\
\p_X & \p_X \ar[l]_{id} \ar[r]^-{(id,\pi)} & \p_X \times_X X}
\end{equation}

From Lemma \ref{lemma_push_forward_fiber_product} and \eqref{eq_gelfand_functorial} we deduce that
\begin{align*}
(\bar p \circ \bar \iota_1)_* \bar \iota_1^* \bar \Phi^*(q_1^* \omega_1 \wedge q_2^* D\omega_2) & = (id,\pi)^* (id \times \pi)_* q_1^*\omega_1 \wedge q_2^*D\omega_2\\
& =\omega_1 \wedge \pi^*\pi_*D\omega_2.
\end{align*}

Similarly,
\begin{displaymath}
\xymatrixcolsep{3pc}
\xymatrix{\p_X \ar[d]^{id} & \bar{\mathcal{N}}_2 \ar[d]^{\bar p \circ \bar \iota_2} \ar[r]^-{\bar \Phi \circ \bar \iota_2} \ar[l]_-{\bar p \circ \bar \iota_2} & \p_X \times_X \p_X \ar[d]^{\pi \times id}\\
\p_X & \p_X \ar[l]_{id} \ar[r]^-{(\pi,id)} & X \times_X \p_X},
\end{displaymath}
commutes and we get

\begin{align*}
(\bar p \circ \bar \iota_2)_* \bar \iota_2^* \bar \Phi^*(q_1^* \omega_1 \wedge q_2^* D\omega_2) & =(-1)^n (\pi,id)^* (\pi \times id)_* q_1^* \omega_1 \wedge q_2^* D\omega_2\\
& = (-1)^n \pi^* \pi_* \omega_1 \wedge D\omega_2.
\end{align*}

From Lemma \ref{lemma_boundary_terms_push_forward} we get
\begin{align*}
d\omega & = \GT(q_1^*D\omega_1 \wedge q_2^* D\omega_2) - \pi^* \pi_*(\omega_1 \wedge s^*D\omega_2) \\
& \quad +\pi^*\pi_* \omega_1 \wedge D\omega_2+D\omega_1 \wedge \pi^*\pi_*\omega_2,
\end{align*}
which is vertical by Lemma \ref{lemma_verticality}.

Note that \eqref{eq_projections_x} implies that
\begin{equation} \label{eq_transform_base_forms}
\GT(q_1^* \pi^*\phi_1 \wedge q_2^*D\omega_2)=\GT(q_1^* D\omega_1 \wedge q_2^* \pi^* \phi_2)=\GT(q_1^* \pi^*\phi_1 \wedge q_2^* \pi^* \phi_2)=0.
\end{equation}
Indeed, tangent vectors to the $n$-dimensional fibers of $\bar p$ are in the kernel of $d(\pi \circ q_1 \circ \bar \Phi)$; and are mapped to horizontal vectors under $d(q_2 \circ \bar \Phi)$.

Hence
\begin{align*}
D\omega+\pi^*\phi & =\GT(q_1^*D\omega_1 \wedge q_2^* D\omega_2)+ \pi^*\pi_* \omega_1 \wedge (D\omega_2+\pi^*\phi_2) + (D\omega_1+\pi^*\phi_1) \wedge \pi^*\pi_*\omega_2 \nonumber \\
& =\GT(q_1^*(D\omega_1+\pi^*\phi_1) \wedge q_2^* (D\omega_2+\pi^*\phi_2)) \nonumber\\
& \quad +\pi^*\pi_* \omega_1 \wedge (D\omega_2+\pi^*\phi_2) + \pi^*\pi_*\omega_2 \wedge (D\omega_1+\pi^*\phi_1).
\end{align*}

The degree of $\GT(q_1^* \omega_1 \wedge q_2^* D\omega_2)$ equals the dimension of the fiber of the map $\pi
\circ \bar p$. Since all vectors which are tangent to the fiber are mapped to horizontal vectors under $q_2 \circ \bar
\Phi$ and since $D\omega_2$ is vertical, it follows that $\pi_*
\GT(q_1^* \omega_1 \wedge q_2^* D\omega_2)=0$. Hence
\eqref{eq_push_forward_omega} follows from the second part of
\eqref{eq_forms_prod} and the projection formula \eqref{projection_formula}.
\endproof

\proof[Proof of Theorem \ref{mthm_prod}]
Let $X$ be an $n$-dimensional (oriented) manifold; let $\mu_i \in \mathcal{V}^\infty(X)$ be represented by $(\omega_i,\phi_i) \in \Omega^{n-1}(\p_X) \times \Omega^n(X)$. We have to show that the product $\mu_1 \cdot \mu_2$ is represented by the forms $\omega,\phi$ in \eqref{eq_forms_prod}.

Let $P \in \mathcal{P}(X)$. Triangulating it, we may assume that $P$
is diffeomorphic to a simplex. We use the forms
$\gamma_0,\kappa,\gamma_1$ which were defined in
\eqref{eq_def_gamma0}-\eqref{eq_def_gamma_i}. By definition of the
product and Theorem \ref{mthm_ext_prod} (and using $[[\Delta
P]]=0$), we obtain

\begin{equation}
\mu_1 \cdot \mu_2(P)=\mu_1 \boxtimes \mu_2(\Delta P)=\int_{\hat N(\Delta P)} (\gamma_0+(-1)^n\kappa)+(-1)^n \int_{N_1(\Delta P)} \gamma_1. \label{eq_reduce_ext_prod}
\end{equation}

From \eqref{eq_transform_base_forms} and the definition of $\GT$ we obtain
\begin{align}
\int_{\bar L_* \bar p^* N(P)} \gamma_0 & = \int_{N(P)} GT(q_1^* \omega_1 \wedge q_2^*(D\omega_2+\pi^* \phi_2)) \nonumber \\
& = \int_{N(P)} GT(q_1^* \omega_1 \wedge q_2^* D\omega_2); \label{eq_main_term}\\
\int_{\bar L_* \bar p^* N(P)} \kappa & = \int_{N(P)} GT(q_1^* \pi^* \phi_1 \wedge q_2^*\omega_2)=0. \label{eq_kappa_term}
\end{align}

By \eqref{eq_diag_r0} we get
\begin{align}
\int_{\rho_* \pi^*P}(\gamma_0+(-1)^n \kappa) & = \int_{\pi^{-1}P} (\omega_1 \wedge s^* (D\omega_2+\pi^* \phi_2)+(-1)^n \pi^*\phi_1 \wedge s^* \omega_2) \nonumber \\
& = \int_P \phi. \label{eq_phi_term}
\end{align}

From \eqref{eq_gelfand_functorial}, \eqref{eq_double_n1}, and Lemmas \ref{lemma_push_forward_fiber_product} and \ref{lemma_normal_cycle_diagonal} we deduce that
\begin{align}
\int_{N_1(\Delta P)} \gamma_1 & = (-1)^{n} (\bar L \circ \bar \iota_1)_* (\bar p \circ \iota_1)^* N(P) (\tau_1^*(q_1^*\omega_1 \wedge q_2^*\omega_2)) \nonumber \\
& = (-1)^{n} \int_{N(P)} (id,\pi)^* (id \times \pi)_* (q_1^*\omega_1 \wedge q_2^*\omega_2) \nonumber \\
& = (-1)^{n} \int_{N(P)} \omega_1 \wedge \pi^*\pi_*\omega_2. \label{eq_gamma1_term}
\end{align}

From \eqref{eq_reduce_ext_prod}, \eqref{eq_main_term}, \eqref{eq_kappa_term}, \eqref{eq_phi_term}, \eqref{eq_gamma1_term} and Lemma \ref{lemma_normal_cycle_diagonal} we get
\begin{align*}
\mu_1 \cdot \mu_2(P)& =\int_{N(P)} \GT(q_1^* \omega_1 \wedge q_2^* D\omega_2)+ \int_{N(P)} \omega_1 \wedge \pi^* \pi_* \omega_2+\int_P \phi,\\
& = \int_{N(P)} \omega+\int_P \phi.
\end{align*}
which finishes the proof of Theorem \ref{mthm_prod}.
\endproof

\section{Functional calculus on valuations}
\label{sec_functional_calc}

We will use Theorem \ref{mthm_prod} in order to introduce a system
of seminorms on $\mathcal{V}^{\infty}(X)$ which behaves well under
the product and defines the usual Fr\'echet topology on
$\mathcal{V}^\infty(X)$ introduced in \cite{ale05d}, Section 3.2.
This will enable us to show Theorem \ref{mthm_functional_calc}.

\proof[Proof of Theorem \ref{mthm_functional_calc}]

Let $\mu$ be a smooth valuation on an $n$-dimensional manifold $X$. Let
$f(z)=\sum_{k=0}^\infty a_k z^k, a_k \in \mathbb{C}$ be an entire function. We have to show that
\begin{equation} \label{eq_power_series}
 f(\mu):=\sum_{k=0}^\infty a_k \mu^k
\end{equation}
converges in $\mathcal{V}^\infty(X)$. Given a compact subset $K
\subset X$, a Riemannian metric $g$ on $X$ and $m \in \mathbb{N}$,
we let $\|\cdot \|_{K,g,m}$ denote the usual $C^m$-seminorm on
smooth forms on $X$ (with respect to $K$ and $g$) and let
$\|\cdot\|_{\pi^{-1}(K),\tilde g,m}$ be the $C^m$-seminorm on forms
on $\p_X$ (with respect to $\pi^{-1}(K)$ and the Sasaki metric
$\tilde g$ on $\p_X$).

We define a semi-norm on $\mathcal{V}^\infty(X)$ by setting
\begin{displaymath}
\|\mu\|_{K,g,m}:=\inf \left\{\|\omega\|_{\pi^{-1}(K),\tilde g,m}+\|\phi\|_{K,g,m}\right\},
\end{displaymath}
where the infimum is over all pairs $(\omega,\phi) \in  \Omega^{n-1}(\p_X) \times \Omega^n(X)$ which represent $\mu$.

This topology is precisely the quotient space topology of
$\Omega^{n-1}(\p_X)\times \Omega^n(X)$. The system of seminorms $\|
\cdot \|_{K,g,m}$ defines the Fr\'echet topology on
$\mathcal{V}^\infty(X)$. Given $(\omega_i,\phi_i) \in
\Omega^{n-1}(\p_X) \times \Omega^n(X)$ ($i=1,2$), define $\omega$
and $\phi$ by \eqref{eq_forms_prod} and set $\mu:=\mu_1 \cdot
\mu_2$. Since the Rumin operator $D$ is a second order differential
operator, we obtain
\begin{displaymath}
\|\omega\|_{\pi^{-1}(K),\tilde g,m} \leq C_{K,g,m} \|\omega_1\|_{\pi^{-1}(K),\tilde g,m} \cdot \|\omega_2\|_{\pi^{-1}(K),\tilde g,m+2}
\end{displaymath}
with a constant $C_{K,g,m}$ depending on $K,g,m$.
Similarly,
\begin{multline*}
\|\phi\|_{K,g,m} \leq C'_{K,g,m} \|\omega_1\|_{\pi^{-1}(K),\tilde g,m} \left(\|\omega_2\|_{\pi^{-1}(K),\tilde g,m+2}+\|\phi_2\|_{K,g,m}\right)\\
+\|\phi_1\|_{K,g,m} \cdot \|\omega_2\|_{\pi^{-1}(K),\tilde g,m}.
\end{multline*}

These inequalities imply
\begin{displaymath}
\|\mu_1\cdot \mu_2\|_{K,g,m} \leq C''_{K,g,m} \|\mu_1\|_{K,g,m}
\cdot \|\mu_2\|_{K,g,m+2}.
\end{displaymath}
In particular, for every smooth valuation $\mu \in \mathcal{V}^\infty(X)$
\begin{displaymath}
\|\mu^2\|_{K,g,m} \leq C''_{K,g,m} \|\mu\|_{K,g,m} \cdot \|\mu\|_{K,g,m+2}.
\end{displaymath}
By iteration
\begin{equation} \label{eq_bound_muk}
\|\mu^k\|_{K,g,m} \leq (C''_{K,g,m})^{k-1} \|\mu\|_{K,g,m} \|\mu\|_{K,g,m+2}^{k-1}.
\end{equation}
We multiply \eqref{eq_bound_muk} by $|a_k|$ and sum over all $k$ to obtain that \eqref{eq_power_series} converges in the
$\| \cdot \|_{K,g,m}$-norm. Since this holds for all $K,g,m$, the theorem follows.
\endproof


\section{Product of generalized valuations}
\label{section_prod_gen}

A generalized valuation on a manifold $X$ is an element of the dual space
\begin{displaymath}
\mathcal{V}^{-\infty}(X):=(\mathcal{V}_c^\infty(X))^*.
\end{displaymath}

Recall that compactly supported valuations can be represented by classes of pairs
$(\omega,\phi) \in \Omega_c^{n-1}(\p_X) \times \Omega_c^n(X)$ modulo the equivalence relation
\begin{equation} \label{eq_kernel_thm}
(\omega_1,\phi_1) \simeq (\omega_2,\phi_2) \iff D\omega_1+\pi^*\phi_1=D\omega_2+\pi^*\phi_2, \pi_*\omega_1=\pi_*\omega_2.
\end{equation}

The embedding $\Xi_\infty:\mathcal{V}^\infty(X) \hookrightarrow \mathcal{V}^{-\infty}(X)$ which is given by
\begin{displaymath}
\Xi_\infty(\mu)(\nu)=\int \mu \cdot \nu
\end{displaymath}
can be described as follows. If $\mu$ is represented by the pair $(\omega,\phi)$, then $\Xi_\infty(\mu)$ is given by the pair of currents
\begin{align*}
T & =\p_X \llcorner s^* (D\omega+\pi^*\phi) \in \mathcal{D}_{n-1}(\p_X),\\
C & = X \llcorner \pi_*\omega \in \mathcal{D}_n(X).
\end{align*}
This follows at once from Theorem \eqref{mthm_prod} (and was earlier proved in \cite{be07b}).

By \eqref{eq_kernel_thm} we obtain that generalized valuations are in one-to-one
correspondence with pairs $(T,C)\in \mathcal{D}_{n-1}(\p_X)\times
\mathcal{D}_n(X)$ satisfying the following three conditions:
\begin{eqnarray}\label{E:pairs1}
T \mbox{ is a cycle, i.e. } \partial T=0;\\
\label{E:pairs2} T \mbox{ is Legendrian};\\\label{E:pairs3}
\pi_*T=\partial C.
\end{eqnarray}

\begin{Lemma}\label{L:generalized-is-smooth}
A generalized valuation $\xi \in \mathcal{V}^{-\infty}(\R^n)$ which corresponds to a pair of $C^\infty$-smooth currents $(T,C)$ belongs to $\mathcal{V}^\infty(\R^n)$.
\end{Lemma}

\proof
Since $T$ is smooth, there is a smooth $n$-form $\kappa$ on $\p_{\R^n}$ such that
\begin{displaymath}
T(\omega)=\int_{\p_{\R^n}} \kappa \wedge \omega
\end{displaymath}
for each $\omega \in \Omega^{n-1}_c(\p_{\R^n})$.

Similarly, since $C$ is smooth, there is some smooth function $\iota$ on $\R^n$ such that
\begin{displaymath}
C(\phi)=\int_{\R^n} \iota \wedge \phi
\end{displaymath}
for each $\phi \in \Omega^{n}_c(\R^n)$.

By conditions \eqref{E:pairs1}, \eqref{E:pairs2} and \eqref{E:pairs3},
$\kappa$ is closed and vertical, and $\pi_* \kappa=(-1)^n d\iota$. Since $s$ is a contactomorphism, $s^* \kappa$ is also closed and vertical. The de Rham cohomology group $H^n(\p_{\R^n})$ vanishes and hence there is some $\omega' \in \Omega^{n-1}(\p_{\R^n})$ with $D\omega'=d\omega'=s^*\kappa$.

Then
\begin{displaymath}
d\pi_*\omega'=\pi_*d\omega'=\pi_*s^*\kappa=(-1)^n \pi_* \kappa=d\iota.
\end{displaymath}

Hence $\pi_*\omega'-\iota \equiv c$ for some constant $c$.

Let $\rho \in \Omega^{n-1}(\p_{\R^n})$ be a transgression form, i.e. $d\rho=0$ and $\pi_*\rho \equiv 1$. Set $\omega:=\omega'-c\rho$. Then $D\omega=D\omega'=s^*\kappa$ and $\pi_*\omega=\iota$. It follows that
$T=\p_{\R^n} \llcorner s^*D\omega$ and $C=X \llcorner \pi_*\omega$ which means, by Theorem \ref{mthm_prod}, that the generalized valuation $\xi$ is represented by $(\omega,0) \in \Omega^{n-1}(\p_{\R^n}) \times \Omega^n(\R^n)$, in particular it is smooth.
\endproof

For closed conic subsets $\Lambda\subset
T^*X\backslash\underline{0},\, \Gamma\subset
T^*\p_X\backslash\underline{0}$, such that
$d\pi_*\Gamma\subset\Lambda$, let us denote by
$\mathcal{V}^{-\infty}_{\Lambda,\Gamma}(X)$ the space of generalized
valuations corresponding to pairs $(T,C)$ such that
\begin{displaymath}
\WF(T)\subset \Gamma,\, \WF(C)\subset \Lambda.
\end{displaymath}

\def\vlg{\mathcal{V}^{-\infty}_{\Lambda,\Gamma}}
Then $\mathcal{V}^{-\infty}_{\Lambda,\Gamma}(X)$ is identified with
a subspace of $C^{-\infty}_\Gamma(\p_X,\Omega^{n-1})\oplus
C^{-\infty}_\Lambda(X)$ which is readily seen to be closed. We equip
$\vlg(X)$ with the induced topology.

\begin{Lemma}\label{C:approxim}
Let $\xi\in \mathcal{V}^{-\infty}_{\Lambda,\Gamma}(\R^n)$ be a
generalized valuation on $\R^n$ given by a pair of currents
$(T,C)\in \mathcal{D}_{n-1}(\p_{\R^n})\times \mathcal{D}_n(\R^n)$,
namely
\begin{displaymath}
WF(C)\subset \Lambda,\, WF(T)\subset \Gamma.
\end{displaymath}

Then there exists a sequence of smooth valuations $\{\xi_j\}\subset
\mathcal{V}^\infty(\R^n)$ corresponding to currents $(T_j,C_j)$ such
that $\xi_j\to\xi$ in $\mathcal{V}^{-\infty}_{\Lambda,\Gamma}(X)$.
\end{Lemma}

\proof

 Let $G$ be the group of all affine transformations of $\R^n$.
Let $\mu_j$ be a sequence of smooth compactly supported measures on
$G$ with integral 1 and whose support converge to the identity
element $e\in G$. Define
\begin{displaymath}
\xi_j:=\int_G(g^*\xi)d\mu_j(g).
\end{displaymath}
Then $\xi_j$ corresponds to the
pair of currents
\begin{displaymath}
(T_j,C_j)=\left(\int_G(g^*T)d\mu_j(g),\int_G(g^*C)d\mu_j(g)\right).
\end{displaymath}

Since $G$ acts transitively on $\R^n$ and on $\p_{\R^n}$ the
currents $C_j$ and $T_j$ are smooth by Proposition
\ref{P:prop-on-wave-fronts}.  Then by Lemma
\ref{L:generalized-is-smooth} $\xi_j\in \mathcal{V}^\infty(\R^n)$.
The rest of the statements follows from Proposition
\ref{P:prop-on-wave-fronts}.

\endproof

The next theorem defines the partial product on generalized
valuations under appropriate technical conditions. These conditions
are formulated in the language when $X$ is affine; however the
conditions do not depend on on a choice of local coordinates, and it
is not hard to rewrite them in an invariant way.

\begin{Theorem} \label{thm_prod_gen}
Let $\Lambda_1,\Lambda_2\subset T^*X\backslash\underline{0},\,
\Gamma_1,\Gamma_2\subset T^*\p_X\backslash\underline{0}$ be closed
conic subsets such that $d\pi_*(\Gamma_i)\subset \Lambda_i,\,
i=1,2$. Let us assume that these subsets satisfy the following
conditions:

\begin{enumerate}
\item[(a)] $\Lambda_1\cap \Lambda_2^s=\emptyset$. \label{condition_down}
\item[(b)] $\Gamma_1\cap (d\pi^*(\Lambda_2))^s=\emptyset$.
\item[(c)] $\Gamma_2\cap (d\pi^*(\Lambda_1))^s=\emptyset$.
\item[(d)] If $(x,[\xi_i],u_i,0) \in
\Gamma_i$ for $i=1,2$, then $u_1 \ne -u_2$.
\item[(e)] If $(x,[\xi])\in \p_X$ and
\begin{align*}
(u,\eta_1) & \in \Gamma_1|_{(x,[\xi])}\\
(-u,\eta_2) & \in \Gamma_2|_{(x,[-\xi])}
\end{align*}
then
\begin{displaymath}
d\theta^*(0,\eta_1,\eta_2) \neq (0,l,-l)\in T^*_{(x,[\xi],[\xi])}(\p_X \times_X \p_X),
\end{displaymath}
where $\theta\colon\p_X \times_X \p_X \to \p_X \times_X \p_X$ is
defined by $\theta(x,[\xi_1],[\xi_2])=(x,[\xi_1],[-\xi_2])$.
\end{enumerate}
Then there is a unique jointly sequentially continuous bilinear map,
called a partial product,
\begin{displaymath}
\mathcal{V}^{-\infty}_{\Lambda_1,\Gamma_1}(X)\times
\mathcal{V}^{-\infty}_{\Lambda_2,\Gamma_2}(X)\to
\mathcal{V}^{-\infty}(X)
\end{displaymath}
satisfying the properties (1)-(3) below.

\begin{enumerate}
\item Commutativity: If $\mu_1,\mu_2 \in \mathcal{V}^{-\infty}(X)$, then
\begin{equation} \label{item_commutativity}
\mu_1 \cdot \mu_2 = \mu_2 \cdot \mu_1
\end{equation}
whenever both sides are defined.
\item Associativity: If $\mu_1,\mu_2, \mu_3 \in \mathcal{V}^{-\infty}(X)$, then
\begin{equation} \label{item_associativity}
(\mu_1 \cdot \mu_2) \cdot \mu_3 = \mu_1 \cdot (\mu_2 \cdot \mu_3)
\end{equation}
whenever both sides are defined.
\item Extension of the smooth product: If $\mu_1,\mu_2 \in \mathcal{V}^\infty(X)$, then $\Xi_\infty(\mu_1) \cdot \Xi_\infty(\mu_2)$ exists and
\begin{equation} \label{eq_extension_smooth}
\Xi_\infty(\mu_1) \cdot \Xi_\infty(\mu_2)=\Xi_\infty(\mu_1 \cdot \mu_2).
\end{equation}
\end{enumerate}
\end{Theorem}

\proof We keep the same notations as in Section \ref{section_prod}.
Let generalized valuations $\mu_i \in
\mathcal{V}^{-\infty}_{\Lambda_i,\Gamma_i}(X), i=1,2,$ be
represented by pairs $(T_i,C_i) \in \mathcal{D}_{n-1}(\p_X) \times
\mathcal{D}_n(X)$ with $T_i$ a Legendrian cycle and
$\pi_*T_i=\partial C_i$. We define
\begin{align}
T & := (-1)^n \GT(q_1^* T_1 \cap q_2^* T_2) + \pi^* C_1 \cap T_2 + T_1 \cap \pi^* C_2 \in \mathcal{D}_{n-1}(\p_X),\label{eq_genproduct_t}\\
C & := C_1 \cap C_2 \in \mathcal{D}_n(X) \label{eq_genproduct_c}
\end{align}
and let $\mu_1 \cdot \mu_2 \in \mathcal{V}^{-\infty}(X)$ be represented by $(T,C)$. This is motivated by the above remark concerning the embedding $\Xi_\infty$ and Proposition \ref{prop_rumin_omega}.

To prove part (1), we have to show that under the conditions (a)-(e)
the currents $T$ and $C$ exist, that $T$ is a Legendrian cycle and
that $\pi_*T=\partial C$.

Since $q_1,q_2$ are submersions, the pull-backs $q_i^*T_i$ are well-defined currents by Proposition \ref{P:lifting}.

By definition,
\begin{displaymath}
q_1^* T_1\cap q_2^*T_2:=(q_1^* T_1 \boxtimes q_2^* T_2)\cap [\p_X \times_X \p_X]
\end{displaymath}
once the right hand side is defined. Using Remark
\ref{restriction-gen-func} and Proposition
\ref{P:wave_front_properties} (iii), this is the case provided
\begin{eqnarray} \label{E:I1}
\left(\Gamma_1\boxtimes \Gamma_2\right)\cap T_{\p_X \times_X
\p_X}^*(\p_X \times \p_X)=\emptyset
\end{eqnarray}
where we denote $\Gamma_1\boxtimes\Gamma_2=(\Gamma_1\times
\Gamma_2)\cup
(\Gamma_1\times\underline{0})\cup(\underline{0}\times\Gamma_2)$.

The fiber of $T_{\p_X \times_X \p_X}^*(\p_X \times \p_X)$ at an arbitrary point $(x,x,[\xi_1],[\xi_2])\in \p_X \times_X \p_X$ is equal to
\begin{eqnarray}\label{E:I2}
T_{\p_X \times_X \p_X}^*(\p_X \times \p_X)|_{(x,x,[\xi_1],[\xi_2])}=\{(u,0,-u,0)|\,
u\in T^*_xX\}.
\end{eqnarray}
Hence (\ref{E:I1}) is equivalent to condition (d).

Note that, if condition (d) is satisfied, by Proposition
\ref{P:wave_front_properties}(iii) a point in $\WF(q_1^*T_1\cap q_2^*T_2)$ has one of the following three forms:
\begin{equation} \label{eq_wave_front_intersection_t}
(x,[\xi_1],[\xi_2],u_1+u_2,\eta_1,\eta_2), (x,[\xi_1],[\xi_2],u_1,\eta_1,0), (x,[\xi_1],[\xi_2],u_2,0,\eta_2),
\end{equation}
where $(x,[\xi_i],u_i,\eta_i) \in \Gamma_i, i=1,2$.

Let us consider the map $\bar\Phi=\Phi\circ \bar L$ where
\begin{displaymath}
\bar \p \stackrel{\bar L}{\longrightarrow} \hat \p \stackrel{\Phi}{\longrightarrow} \p_X \times_X \p_X.
\end{displaymath}

Recall that $\Phi$ is a submersion with fibers diffeomorphic to a
closed segment, and $\bar L$ is the oriented blow up of $\hat \p$ along
the smooth submanifold (without boundary)
\begin{displaymath}
 A:=\{(x,[\xi:-\xi],[\xi],[-\xi])\} \subset \hat \p.
\end{displaymath}

Let us denote by $\Delta$ the diagonal
\begin{displaymath}
\Delta:=\{(x,[\xi],[\xi])\}\subset \p_X \times_X \p_X.
\end{displaymath}

Fix an arbitrary point $a=(x,[\xi:-\xi],[\xi],[-\xi])\in A$ and set
\begin{align*}
b & :=\Phi(a)=(x,[\xi],[-\xi])\in \p_X \times_X \p_X,\\
c & :=\theta(b)=(x,[\xi],[\xi]) \in \Delta.
\end{align*}

Clearly the $\Phi$-image of $A$ is a smooth submanifold $\Phi(A)$
which is equal to $\theta(\Delta)$ (where, we recall, $\theta$ is
the involution $\theta(x,[\xi_1],[\xi_2])=(x,[\xi_1],[-\xi_2])$).
Thus
\begin{eqnarray} \label{E:5}
d\Phi(T_aA)=d\theta(T_{\theta(a)}\Delta).
\end{eqnarray}

Note that
\begin{displaymath}
T^*_\Delta(\p_X \times_X \p_X)|_c=\{(0,l,-l)|\, l\in T_{[\xi]}^*\p_X\}.
\end{displaymath}

If $\zeta \in \WF(q_1^*T_1 \cap q_2^*T_2)|_b$ is of the form $(x,[\xi],[-\xi],u_1,\eta_1,0)$ or $(x,[\xi],[-\xi],u_2,0,\eta_2)$, then clearly $d\theta^*(\zeta)\not\in T^*_\Delta(\p_X \times_X \p_X)|_c$.

If $\zeta$ is of the form $(x,[\xi],[-\xi],u_1+u_2,\eta_1,\eta_2)$
with $(x,[\xi],u_1,\eta_1) \in \Gamma_1$ and $(x,[-\xi],u_2,\eta_2)
\in \Gamma_2$, then condition (e) of Theorem \ref{thm_prod_gen}
implies that $d\theta^*(\zeta)\not\in T^*_\Delta(\p_X \times_X
\p_X)|_c$.

Using \eqref{eq_wave_front_intersection_t} we get that for any $\zeta\in \WF(q_1^*T_1 \cap q_2^*T_2)|_b$ the
restriction of $\zeta$ to the subspace (\ref{E:5}) does not vanish. Thus by Proposition \ref{prop_P:2} the pullback
$\bar\Phi^*(q_1^*T_1 \cap q_2^*T_2)$ is well defined, which implies that the first term in \eqref{eq_genproduct_t} is defined.

By Propositions \ref{P:lifting}, \ref{P:tens-prod-wf} and conditions (b) and (c), the intersections $\pi^* C_1 \cap T_2$ and $T_1 \cap \pi^* C_2$ exist. Similarly, condition (a) implies that the intersection $C_1 \cap C_2$ exists.

This shows that $T$ and $C$ are well-defined. The joint sequential
continuity of the map
\begin{displaymath}
 ((T_1,C_1),(T_2,C_2))\mapsto (T,C)
\end{displaymath}
follows from the sequential continuity properties in Propositions
\ref{P:lifting}, \ref{P:tens-prod-wf}, and \ref{P:wf-push-forward}.

Before we show that $T$ is a Legendrian cycle and $\pi_*T=C$, let us
prove \eqref{eq_extension_smooth}. Suppose that $\mu_i, i=1,2$ is
presented by $(\omega_i,\phi_i) \in \Omega^{n-1}(\p_X) \times
\Omega^n(X)$ and that the product $\mu_1 \cdot \mu_2$ is represented
by $(\omega,\phi)$. Setting
\begin{align*}
T_i & :=[\p_X] \llcorner s^*(D\omega_i+\phi_i)\\
T & := [\p_X] \llcorner s^*(D\omega+\phi)\\
C_i & :=[X] \llcorner \pi_*\omega_i\\
C & := [X] \llcorner \pi_*\omega
\end{align*}
we have to show that \eqref{eq_genproduct_t} and
\eqref{eq_genproduct_c} are satisfied.

Since $\pi_*\omega=\pi_*\omega_1 \wedge \pi_* \omega_2$ by
\eqref{eq_push_forward_omega}, we clearly get $C=C_1 \cap C_2$. Now
we compute
\begin{align*}
\GT(q_1^*T_1 \cap q_2^*T_2) & = \bar p_* \bar \Phi^* ([\p_X \times_X \p_X] \llcorner q_1^*s^*(D \omega_1+\pi^*\phi_1) \wedge q_2^*s^*(D \omega_2+\pi^*\phi_2)) \\
& = \bar p_* [\bar \p] \llcorner \bar \Phi^* (q_1^*s^*(D \omega_1+\pi^*\phi_1) \wedge q_2^*s^*(D \omega_2+\pi^*\phi_2))\\
& = (-1)^{n(3n-1-2n)}[\p_X] \llcorner \bar p_* \bar \Phi^* (q_1^*s^*(D \omega_1+\pi^*\phi_1) \wedge q_2^*s^*(D \omega_2+\pi^*\phi_2))\\
& = [\p_X] \llcorner \GT(q_1^*s^*(D \omega_1+\pi^*\phi_1) \wedge q_2^*s^*(D \omega_2+\pi^*\phi_2))\\
& = (-1)^n [\p_X] \llcorner s^* \GT(q_1^*(D \omega_1+\pi^*\phi_1)
\wedge q_2^*(D \omega_2+\pi^*\phi_2)).
\end{align*}

Next,
\begin{align*}
\pi^*C_1 \cap T_2 & = [\p_X] \llcorner \pi^*\pi_* \omega_1 \cap [\p_X] \llcorner s^*(D\omega_2+\pi^*\phi_2)\\
& = [\p_X] \llcorner s^*(\pi^*\pi_*\omega_1 \wedge
(D\omega_2+\pi^*\phi_2))
\end{align*}
and similarly
\begin{displaymath}
T_1 \cap \pi^*C_2 = [\p_X] \llcorner s^*(\pi^*\pi_*\omega_2 \wedge
(D\omega_1+\pi^*\phi_1)).
\end{displaymath}

These three equations and  \eqref{eq_Rumin_omega} imply
\begin{displaymath}
(-1)^n \GT(q_1^*T_1 \cap q_2^*T_2)+\pi^*C_1 \cap T_2+T_1 \cap
\pi^*C_2=[\p_X] \llcorner s^*(D\omega+\pi^*\phi)=T.
\end{displaymath}

Now let us return to the situation of generalized valuations. Recall
that we have to show that if $\mu_i\in
\mathcal{V}^{-\infty}_{\Lambda_i,\Gamma_i}(X)$, $i=1,2$, and $T,C$
are defined by \eqref{eq_genproduct_t}, \eqref{eq_genproduct_c} then
$T$ is a Legendrian cycle and $\pi_*T=C$. First observe that if
$\mu_i$ are smooth valuations then this is true because of just
proven compatibility of the partial product on generalized
valuations with the product on smooth valuations (in this case the
pair $(T,C)$ corresponds to the valuation $\mu_1\cdot\mu_2$). In
general, the question is local hence we may and will assume that
$X=\R^n$. By Lemma \ref{C:approxim} there exist two sequences of
smooth valuations $\mu_{i,j}\to \mu_i$ in
$\mathcal{V}^{-\infty}_{\Lambda_i,\Gamma_i}(X)$ as $j\to\infty$,
$i=1,2$. Then by joint sequential continuity $\mu_{1,j}\cdot
\mu_{2,j}\to \mu_1\cdot\mu_2$ in $\mathcal{V}^{-\infty}(X)$. This
means that when $j\to \infty$ one has $C_{i,j}\to C_i$ in
$\mathcal{D}_{n-1}(X)$, $T_{i,j}\to T_i$ in $\mathcal{D}_{n}(\p_X)$.
The result follows by continuity.

Commutativity and associativity of the product on generalized
functions (i.e. equations \eqref{item_commutativity} and
\eqref{item_associativity}) are easily proven, using the fact that
the intersection product $\cap$ is (graded) commutative and
associative whenever it is defined.

The uniqueness of the partial product follows from joint sequential
continuity of it, property (3) of the theorem, and the fact that
$\mathcal{V}^\infty(X)$ is sequentially dense in
$\mathcal{V}^\infty_{\Lambda,\Gamma}(X)$ by Lemma \ref{C:approxim}.

\endproof

\section{Product and intersection of manifolds with corners}
\label{sec_prod_man_corners}

In this section, we give a precise form to the statement that the product of valuations corresponds to taking the intersection of sets.
In a slightly more general form, this was conjectured in \cite{ale06}.

The space $\mathcal{P}(X)$ of compact submanifolds with corners is
also embedded in $\mathcal{V}^{-\infty}(X)$, by the map
\begin{align*}
 \Xi_\mathcal{P}: \mathcal{P}(X) & \hookrightarrow \mathcal{V}^{-\infty}(X) \\
P & \mapsto (\mu \mapsto \mu(P)).
\end{align*}

\proof[Proof of Theorem \ref{mthm_mflds_corners}] Let
$P^{(1)},P^{(2)}$ be compact submanifolds with corners which
intersect transversally. We have to show that the product of the
generalized valuations $\Xi_\mathcal{P}(P^{(1)})$ and
$\Xi_\mathcal{P}(P^{(2)})$ exists and equals
$\Xi_\mathcal{P}(P^{(1)} \cap P^{(2)})$.

Let us first show existence.

Let us fix two submanifolds with corners $P^{(1)},P^{(2)} \subset X$.
Assume moreover that they intersect transversally in sense of
Definition \ref{D:transversality}. Let
\begin{align*}
 T_i & :=N(P^{(i)}) \in D_{n-1}(\p_X),\\
 C_i & :=P^{(i)} \in D_n(X), \quad  i=1,2.
\end{align*}

We have to check the conditions (a)-(e) of Theorem \ref{thm_prod_gen}.

Let us fix an arbitrary point $x_0\in P^{(1)} \cap P^{(2)}$.
Since $P^{(i)}$ are submanifolds with corners, in a neighborhood of
$x_0$ there exist smooth functions
\begin{displaymath}
 f_1^{(i)},\dots, f_{k_i}^{(i)},g_1^{(i)},\dots,g_{l_i}^{(i)},\,
i=1,2,
\end{displaymath}
such that
\begin{itemize}
 \item these functions vanish at $x_0$;
\item in a neighborhood of $x_0$
\begin{displaymath}
P^{(i)}=\{x|\, f_1^{(i)}(x)=\dots =f_{k_i}^{(i)}(x)=0,\, g_1^{(i)}(x)\geq
0,\dots, g_{l_i}^{(i)}(x)\geq 0\};
\end{displaymath}
\item for each $i=1,2$ the sequence
\begin{eqnarray} \label{E:seqi}
\nabla f_1^{(i)}(x_0),\dots,\nabla f_{k_i}^{(i)}(x_0), \nabla
g_1^{(i)}(x_0),\dots, \nabla g_{l_i}^{(i)}(x_0)\in T^*_{x_0}X
\end{eqnarray}
is linearly independent.
\end{itemize}

By the assumption of transversality of $P^{(1)},P^{(2)}$ the union of
sequences (\ref{E:seqi}) over $i=1,2$
\begin{eqnarray*}\label{E:sequence}
\nabla f_1^{(1)}(x_0),\dots,\nabla f_{k_1}^{(1)}(x_0), \nabla
g_1^{(1)}(x_0),\dots, \nabla g_{l_1}^{(1)}(x_0),\\\nabla
f_1^{(2)}(x_0),\dots,\nabla f_{k_2}^{(2)}(x_0), \nabla
g_1^{(2)}(x_0),\dots, \nabla g_{l_2}^{(2)}(x_0)
\end{eqnarray*}
is linearly independent. We can choose smooth functions
$h_1,\dots,h_e$ such that the union of the former sequence with
$\nabla h_1(x_0),\dots,\nabla h_e(x_0)$ form a basis of
$T^*_{x_0}X$. Then in a neighborhood of $x_0$ all the functions
\begin{displaymath}
 \{f^{(1)}_j\}, \{f^{(2)}_j\}, \{g^{(1)}_j\},\{g^{(2)}_j\},\{h_j\}
\end{displaymath}
form a coordinate system.

Thus we may assume that
\begin{align}
X & =\R^a\times \R^b \times \R^c \times \R^d \times \R^e,\\
P^{(1)} & =0 \times \R^b \times \R^c_{\geq 0} \times
\R^d \times \R^e,\label{E:Aone} \\
P^{(2)} & = \R^a \times 0 \times \R^c \times \R^d_{\geq 0} \times \R^e \label{E:Atwo}
 \end{align}
and $x_0=0.$ Here we abuse the notation and denote
$a=k_1,b=k_2,c=l_1,d=l_2$.

To check (a) we note that the fibers of the wave front at $0$ satisfy
\begin{eqnarray}
\WF(P^{(1)})|_{0} \subset \R^{a*} \times 0 \times\R^{c*} \times 0 \times
0,\\
\WF(P^{(2)})|_{0}\subset 0 \times \R^{b*} \times 0 \times \R^{d*} \times
0.
\end{eqnarray}
It follows that $\WF(P^{(1)})|_{0}\cap \WF(P^{(2)})^s|_{0}=\emptyset$, which is condition (a).

Let us denote
\begin{align*}
L & :=\R^a\times \R^c,\\
M & :=\R^b\times \R^d\times \R^e.
\end{align*}

Thus $\R^n=L\times M$. By Proposition \ref{prop_P:2}, for any $(x,[\xi])\in \p_X$
\begin{displaymath}
\WF(\pi^*C_1)|_{(x,[\xi])}\subset d\pi^*(\WF(C_1))=\{(l,0)|\, l\in
\WF(C_1)|_x\}.
\end{displaymath}

Hence
\begin{eqnarray}\label{E:middle1}
\WF(\pi^*C_1)|_{(0,[\xi])} \subset (L^* \times 0_{M^*}) \times
0_{T^*_{[\xi]}\p_+(\R^{n*})}
\end{eqnarray}
where $0_{T^*_{[\xi]} \p_+(\R^{n*})}$ denotes the zero vector in
the cotangent space $T^*_{[\xi]}\p_+(\R^{n*})$.

Next $P^{(2)}=L\times M_1$ where $M_1=0 \times \R^{d}_{\geq 0}
\times \R^e$. Hence
\begin{displaymath}
N(P^{(2)})=L \times N(M_1) \subset \p_+(T^*(L\times M)).
\end{displaymath}

Then by Proposition \ref{P:wave_front_properties}(iii)
\begin{eqnarray}\label{E:middle2}
\WF(T_2)|_{(0,[\xi])} \subset (0_{L^*} \times M^*) \times
T^*_{[\xi]}(\p_+(\R^{n*})).
\end{eqnarray}
Clearly (\ref{E:middle1})-(\ref{E:middle2}) imply that
\begin{displaymath}
 \WF(\pi^*C_1)|_{(0,[\xi])}\cap
\WF(T_2)^s|_{(0,[\xi])}=\emptyset.
\end{displaymath}

This is condition (c); condition (b) follows by symmetry.

To check the condition (d) let us assume that
$(0,[\xi_1],[\xi_2]) \in \p_X$ and
\begin{eqnarray*}
(u,0)\in \WF(T_1)|_{(0,[\xi_1])} \subset T^*_{(0,[\xi_1])}(\p_X),\\
(-u,0)\in \WF(T_2)|_{(0,[\xi_2])} \subset T^*_{(0,[\xi_2])}(\p_X).
\end{eqnarray*}

Using Proposition \ref{P:wave_front_properties}(iii), the equations
\eqref{E:Aone}, \eqref{E:Atwo} imply respectively
\begin{displaymath}
u\in \R^{a*} \times 0 \times \R^{c*} \times 0\times 0,\\\label{E:u2}
-u\in 0\times\R^{b*}\times 0\times \R^{d*}\times 0
\end{displaymath}
and hence $u=0$. But this is
impossible since $(u,0)\in \WF(T_1)$. Thus condition
(d) is satisfied.

Let us check conditions (e). Let us fix
$(0,[\xi])\in \p_X$. Assume that
\begin{eqnarray*}
(u,\eta_1) \in \WF(T_1)|_{(0,[\xi])},\\
(-u,\eta_2)\in \WF(T_2)|_{(0,[\xi])}.
\end{eqnarray*}

For an $\R_{>0}$-invariant set $U$ let us denote
$\p_+(U):=(U\backslash \{0\})/\R_{>0}$. It is easy to see that the
fibers of $N(P^{(i)})$ over $0$ are equal to
\begin{align*}
N(P^{(1)})|_{0} & = \p_+ \left(\R^{a*} \times 0 \times \R^{c*}_{\geq
0} \times 0 \times 0 \right),\\
N(P^{(2)})|_{0} & =\p_+\left(0 \times \R^{b*} \times
0 \times\R^{d*}_{\geq 0} \times 0\right).
\end{align*}

We easily conclude that
\begin{displaymath}
[\xi]\in \p_+\left(\R^{a*} \times 0 \times \R^{c*}_{\geq
0} \times 0 \times 0\right) \cap \p_+\left(0 \times \R^{b*} \times
0 \times \R^{d*}_{\geq 0} \times 0\right)=\emptyset.
\end{displaymath}
This is a contradiction. Hence condition (e) is vacuous.

This establishes the existence of the product
$\Xi_\mathcal{P}(P^{(1)}) \cdot \Xi_\mathcal{P}(P^{(2)})$. It
remains to show that this product equals $\Xi_\mathcal{P}(P^{(1)}
\cap P^{(2)})$.

We define $T$ and $C$ by \eqref{eq_genproduct_t} and
\eqref{eq_genproduct_c} and want to see that $T=N(P^{(1)} \cap
P^{(2)}), C=P^{(1)} \cap P^{(2)}$. The second equation is obvious.
For the first one, we let $(x,[\zeta]) \in \spt \GT(q_1^*T_1 \cap
q_2^*T_2)$. Then there exist $(x,[\xi]) \in \spt T_1$ and
$(x,[\eta]) \in \spt T_2$ with $\xi+\eta=\zeta$. This implies
$(x,[\zeta]) \in N(P_1 \cap P_2)$.

If $(x,[\zeta]) \in \spt \pi^*C_1 \cap T_2$, then $x \in P_1$ and
$(x,[\zeta]) \in \spt T_2$. Again we deduce that $(x,[\zeta]) \in
N(P^{(1)} \cap P^{(2)})$. A similar argument gives $\spt T_1 \cap
\pi^*C_2 \subset N(P^{(1)} \cap P^{(2)})$.

We have shown that
\begin{displaymath}
\spt T \subset N(P^{(1)} \cap P^{(2)}).
\end{displaymath}

By Lemma \ref{lemma_constancy_theorem},

\begin{displaymath}
T=m N(P^{(1)} \cap P^{(2)})
\end{displaymath}
for some integer $m$.

In order to fix $m$, we compute
\begin{displaymath}
m \pi_* N(P^{(1)} \cap P^{(2)})=\pi_*T=\partial C=\partial (P^{(1)}
\cap P^{(2)}) = \pi_* N(P^{(1)} \cap P^{(2)}).
\end{displaymath}
Hence $m=1$.
\endproof

\section{Kinematic formulas in compact rank one symmetric spaces}
\label{section_kin_form}

Let $X=G/H$ be a compact rank one symmetric space. Since $X$ is two-point homogeneous, the induced action of $G$ on $\p_X=\p_+(T^*X)$ is transitive.

Since we can average over the $G$-action, every smooth $G$-invariant valuation $\mu$ on $X$ is represented by a pair $(\omega,\phi) \in \Omega^{n-1}(\p_X) \times \Omega^n(X)$ of $G$-invariant forms. In particular, $\mathcal{V}^\infty(X)$ is finite-dimensional. Let $\varphi_1,\ldots,\varphi_N$ be a basis of $\mathcal{V}^\infty(X)^G$.

Let $dg$ be the unique Haar measure of volume $1$ on $G$. Fu, using a slightly different language, showed in \cite{fu90} that for each $\mu \in \mathcal{V}^\infty(X)^G$, there exist constants $c_{ij}^\mu \in \mathbb{C}$ such that for all $P_1,P_2 \in\mathcal{P}(X)$ the following kinematic formula holds true:
\begin{equation} \label{eq_kinematic_formula}
\int_G \mu(P_1 \cap gP_2)dg=\sum_{i,j=1}^N c_{ij}^\mu \varphi_i(P_1)\varphi_j(P_2).
\end{equation}

Similarly as in \cite{fu06}, we define a map
\begin{align*}
k_G:\mathcal{V}^\infty(X)^G & \to \mathcal{V}^\infty(X)^G \otimes \mathcal{V}^\infty(X)^G,\\
\mu & \mapsto \sum_{i,j=1}^N c_{ij}^\mu \varphi_i \otimes \varphi_j.
\end{align*}

Let
\begin{displaymath}
m_G: \mathcal{V}^\infty(X)^G \otimes \mathcal{V}^\infty(X)^G \to \mathcal{V}^\infty(X)^G
\end{displaymath}
be the restriction of the multiplication to smooth $G$-invariant valuations.

The Poincar\'e duality proved in \cite{ale05d} implies that there is an isomorphism
\begin{displaymath}
p: \mathcal{V}^\infty(X) \to \mathcal{V}^{-\infty}(X)
\end{displaymath}
such that
\begin{displaymath}
\langle p(\mu),\nu\rangle=\int \mu \cdot \nu=(\mu \cdot \nu)(X), \quad \mu,\nu \in \mathcal{V}^\infty(X).
\end{displaymath}
With $i_G: \mathcal{V}^\infty(X)^G \hookrightarrow \mathcal{V}^\infty(X)$ being the natural inclusion map, we let $p_G:=i_G^* \circ p \circ i_G:\mathcal{V}^\infty(X)^G  \to \mathcal{V}^\infty(X)^{G*}$.

\proof[Proof of Theorem \ref{mthm_kinematic}]

Let us fix a basis $\varphi_1,\ldots,\varphi_N$ of $\mathcal{V}^\infty(X)^G$. We set $g_{ij}:=\langle p_G \varphi_i,\varphi_j\rangle$. This is an invertible symmetric matrix, whose inverse will be denoted by $(g^{ij})$.

For $\phi \in \mathcal{V}^\infty(X)$, the valuation $\phi^G:=\int_G g^*\phi dg$ is invariant. Hence we can express it in terms of our basis as
\begin{displaymath}
 \phi^G=\sum_{i,j=1}^N \langle p_G \varphi_j,\phi^G\rangle g^{ji} \varphi_i=\sum_{i,j=1}^N \langle p \varphi_j,\phi\rangle g^{ji} \varphi_i.
\end{displaymath}

For $\nu \in \mathcal{V}^{-\infty}(X)^G$ we obtain that
\begin{displaymath}
\langle \nu,\phi\rangle = \langle \nu,\phi^G \rangle = \sum_{i,j=1}^N  g^{ji} \langle \nu,  \varphi_i\rangle \langle p\varphi_j,\phi\rangle
\end{displaymath}
and thus
\begin{equation} \label{eq_ginv_implies_smooth}
\nu = p\left(\sum_{i,j=1}^N  g^{ji} \langle \nu,  \varphi_i\rangle \varphi_j\right)
\end{equation}
is smooth. In other words, every $G$-invariant generalized valuation is smooth.

We have to show that the following diagram commutes
\begin{displaymath}
\xymatrixcolsep{3pc}
\xymatrix{\mathcal{V}^\infty(X)^G \ar[d]^{k_G} \ar[r]^{p_G}& \mathcal{V}^\infty(X)^{G*} \ar[d]^{m_G^*}\\
\mathcal{V}^\infty(X)^G \otimes \mathcal{V}^\infty(X)^G \ar[r]^-{p_G \otimes p_G} & \mathcal{V}^\infty(X)^{G*} \otimes \mathcal{V}^\infty(X)^{G*}},
\end{displaymath}

Fix $\mu \in \mathcal{V}^\infty(X)^G$ and $P_1, P_2 \in \mathcal{P}(X)$. Whenever $g_1,g_2$ are such that $g_1P_1 \cap g_2P_2$ intersect transversally strata by strata we can apply Theorem \ref{mthm_mflds_corners} to get
\begin{displaymath}
\mu(g_1P_1 \cap g_2 P_2)= \langle \Xi_\mathcal{P}(g_1P_1 \cap g_2P_2),\mu\rangle=\left\langle \Xi_\mathcal{P}(g_1P_1) \cdot \Xi_\mathcal{P}(g_2P_2),\mu \right\rangle.
\end{displaymath}

Integrating this equality over $G \times G$ we find
\begin{align}
\int_G \mu(P_1 \cap gP_2)dg & = \int_{G \times G} \mu(g_1P_1 \cap g_2 P_2) d(g_1,g_2) \nonumber \\
& = \left \langle \int_G \Xi_\mathcal{P}(g_1P_1)dg_1 \cdot \int_G \Xi_{\mathcal{P}}(g_2P_2)dg_2,\mu\right\rangle. \label{eq_integral_formula}
\end{align}

By \eqref{eq_ginv_implies_smooth} we have
\begin{displaymath}
\int_G \Xi_\mathcal{P}(gP)dg=\sum_{i,j=1}^N g^{ij}\varphi_j(P)p \varphi_i.
\end{displaymath}

Plugging this into \eqref{eq_integral_formula} and using \eqref{eq_extension_smooth}, we find
\begin{align*}
\int_G \mu(P_1 \cap gP_2)dg  & = \left\langle p\left(\sum_{i,k=1}^N  g^{ik} \varphi_i(P_1) \varphi_k \cdot \sum_{j,l=1}^N g^{jl} \varphi_j(P_2) \varphi_l\right),\mu\right\rangle\\
& = \sum_{i,j=1}^N \left(\sum_{k,l=1}^N g^{ik}g^{jl} \langle p(\varphi_k \cdot \varphi_l),\mu\rangle\right) \varphi_i(P_1)\varphi_j(P_2).
\end{align*}

Comparing with \eqref{eq_kinematic_formula}, we see that
\begin{displaymath}
c_{i,j}^\mu=\sum_{k,l=1}^N g^{ik}g^{jl} \langle p(\varphi_k \cdot \varphi_l),\mu\rangle=\sum_{k,l=1}^N g^{ik}g^{jl} \langle p_G \mu,\varphi_k \cdot \varphi_l\rangle.
\end{displaymath}

For $s,t=1,\ldots,N$ we obtain
\begin{align*}
p_G \otimes p_G \circ k_G(\mu)(\varphi_s,\varphi_t) & = \sum_{i,j} c_{i,j}^\mu \langle p_G \varphi_i,\varphi_s\rangle \langle p_G \varphi_j,\varphi_t\rangle\\
& = \sum_{i,j,k,l} g^{ik}g^{jl} g_{is}g_{jt} \langle p_G\mu,\varphi_k\cdot \varphi_l\rangle\\
& = \langle p_G \mu,\varphi_s \cdot \varphi_t\rangle\\
& = m_G^* \circ p_G(\mu)(\varphi_s,\varphi_t\rangle.
\end{align*}
\endproof

Remarks:
\begin{enumerate}
\item A similar theorem in the translation invariant case was proven in \cite{befu06},
based on \cite{fu06}. Roughly speaking, the algebra structure on $\mathcal{V}^\infty(X)^G$
(respectively $\Val^G$ in the translation invariant case) determines the kinematic formulas for the group $G$ and vice versa.
This observation was used in \cite{befu08} to write down explicitly the principal kinematic formula for hermitian vector spaces.
We hope that Theorem \ref{mthm_kinematic} will be the starting point for a
similar treatment of kinematic formulas on compact rank one symmetric spaces.
\item We used the compactness hypothesis at several places. It would be interesting to have an analogue of
Theorem \ref{mthm_kinematic} in the non-compact case.
\end{enumerate}



\begin{thebibliography}{99}
\bibitem{ale01}  Alesker, S.: Description of translation invariant valuations on convex sets with solution of P. McMullen's conjecture.
{\it Geom. Funct. Anal.} \textbf{11} (2001), 244--272.

\bibitem{ale04}  Alesker, S.: The multiplicative structure on polynomial valuations. {\it Geom. Funct. Anal.} \textbf{14} (2004), 1--26.

\bibitem{ale05a} Alesker, S.: Theory of valuations on manifolds, I. Linear Spaces. {\it Israel J. Math.} \textbf{156} (2006), 311--340.

\bibitem{ale05b} Alesker, S.: Theory of valuations on manifolds, II. {\it Adv. Math.} \textbf{207} (2006), 420--454.

\bibitem{ale05d} Alesker, S.: Theory of valuations on manifolds
IV. New properties of the multiplicative structure. In {\it
Geometric aspects of functional analysis}, LNM \textbf{1910} (2007),
1--44.

\bibitem{ale06} Alesker, S.: Valuations on manifolds: a survey. {\it Geom. Funct. Anal.} \textbf{17} (2007), 1321--1341.

\bibitem{ale05c} Alesker, S., Fu, J. H. G.: Theory of valuations on
  manifolds, III. Multiplicative structure in the general case. {\it Trans. Amer. Math. Soc.}  \textbf{360}  (2008), 1951--1981.

\bibitem{alfe98} \'Alvarez Paiva, J. C., Fernandes, E.: Crofton formulas in
  projective Finsler spaces. {\it Electron. Res. Announc. Am. Math. Soc.}
  \textbf{4} (1998), 91--100.

\bibitem{alfe05} \'Alvarez Paiva, J. C., Fernandes, E.: Crofton formulas and
  Gelfand transforms. {\it Selecta Math.} \textbf{13} (2007), 369--390.

\bibitem{begeve92} Berline, N., Getzler, E., Vergne, M.: {\it Heat kernels
  and Dirac operators}. Springer-Verlag Berlin, 1992.

\bibitem{be06} Bernig, A.: Support functions, projections and
  Minkowski addition of Legendrian cycles. {\it Indiana U. Math. J.}
  \textbf{55} (2006), 443--464.

\bibitem{be07} Bernig, A.: Valuations with Crofton formula and Finsler geometry. {\it Adv. Math.} \textbf{210} (2007), 733--753.

\bibitem{be07b} Bernig, A.: A product formula for valuations on manifolds with applications to the integral geometry of the quaternionic line. {\it Comm. Math. Helv.} \textbf{84} (2009), 1--19.

\bibitem{be08a} Bernig, A.: A Hadwiger-type theorem for the special unitary group. {\it Geom. Funct. Anal.} \textbf{19} (2009), 356--372.

\bibitem{be08b} Bernig, A.: Integral geometry under $G_2$ and $Spin(7)$. arXiv:0803.3885. To appear in {\it Israel J. Math.}

\bibitem{bebr07} Bernig, A., Br\"ocker, L.: Valuations on manifolds
  and Rumin cohomology. {\it J. Differential Geom.} \textbf{75} (2007), 433--457.

\bibitem{befu06} Bernig, A., Fu. J. H. G.: Convolution of valuations.  {\it Geom. Dedicata} \textbf{123} (2006), 153--169.

\bibitem{befu08} Bernig, A., Fu, J. H. G.: Hermitian integral geometry.  arXiv:0801.0711. To appear in {\it Ann. of Math.}

\bibitem{bes87} Besse, A. L.: {\it Einstein manifolds.} Ergebnisse der Mathematik und ihrer Grenzgebiete. 3. Folge, Bd. 10. Berlin etc.: Springer-Verlag (1987).

\bibitem{bro66} Brothers, J. E.: Integral geometry in homogeneous
  spaces. {\it Trans. Amer. Math. Soc.} \textbf{124} (1966), 480--517.


\bibitem{fed69} Federer, H.: {\it Geometric Measure Theory}.
Springer-Verlag New York 1969.

\bibitem{fu90}  Fu, J. H. G.:  Kinematic formulas in integral geometry.
{\it Indiana Univ. Math. J.} \textbf{39} (1990), 1115--1154.

\bibitem{fu06} Fu, J. H. G.: Structure of the unitary valuation algebra. {\it J. Differential Geom.} \textbf{72} (2006), 509--533.

\bibitem{gui-st77} Guillemin, V., Sternberg, Sh.: {\it Geometric asymptotics}.
Mathematical Surveys 14. American Mathematical Society (AMS), Providence, R.I. 1977.

\bibitem{helgason}
Helgason, S.: {\it Differential geometry and symmetric spaces}. Pure
and Applied Mathematics, Vol. XII. Academic Press, New York-London
1962.

\bibitem{hormander}
H\"ormander, L.: {\it The analysis of linear partial differential
operators. I. Distribution theory and Fourier analysis}. Grundlehren
der Mathematischen Wissenschaften [Fundamental Principles of
Mathematical Sciences], 256. Springer-Verlag, Berlin, 1983.

\bibitem{rum94} Rumin, M.: Formes diff\'erentielles sur les vari\'et\'es de contact.
{\it J. Differential Geom.} \textbf{39} (1994), 281--330.

\end{thebibliography}
\end{document}